\documentclass[11pt,a4paper]{article}

\usepackage{amsmath,amssymb}
\newcommand{\e}{\varepsilon}
\newcommand{\al}{\alpha}
\newcommand{\D}{\Delta}

\newcommand{\n}{\nabla}

\newcommand{\g}{\int_{\mathbb{R}^{d}}}
\newcommand{\p}{\partial}
\newcommand{\va}{\varphi}

\newcommand{\R}{\mathbb{R}}

\newtheorem{definition}{Definition}
\newtheorem{notation}{Notation}
\newtheorem{theorem}{Theorem}
\newtheorem{proposition}{Proposition}
\newtheorem{corollaire}{Corollary}

\newtheorem{rem}{Remark}
\newtheorem{lemme}{Lemma}

\title{Study of compactness for compressible fluid models with a specific Korteweg tensor}
\author{Boris Haspot \thanks{Universit\'e Paris XII - Val de Marne 61, avenue
du G\'en\'eral de Gaulle 94 010 CRETEIL Cedex T\'el\'ephone : (33-1)
45 17 16 51 T\'el\'ecopie : (33-1) 45 17 16 49 e-mail :
haspot@univ-paris12.fr}}
\date{}
\begin{document}
\maketitle
\subsubsection*{Abstract}
This work is devoted to the global stability theory of solutions for
a general  isothermal model of capillary fluids
derived by C. Rohde in \cite{4Ro}, which can be used as a phase transition model.\\
This chapter is structured in the following way: first of all
inspired by the result by P.-L. Lions in \cite{4L2} on the
Navier-Stokes compressible system we will show the global stability
of weak solution for our system with isentropic pressure and next
with general pressure. Next we will consider perturbations close to
a stable equilibrium as in the case of strong solution.
\section{Introduction}
\subsection{Presentation of the model}
The correct mathematical description of liquid-vapor phase
interfaces and their dynamical behavior in compressible fluid flow
has a long history. We are concerned with compressible fluids
endowed with internal capillarity. One of the first model which
takes into consideration the variation of density on the interface
between two phases, originates from the XIXth century work by Van
der Waals and Korteweg \cite{4Ko}. It was actually derived in his
modern form in the 1980s using the second gradient theory, see for
instance \cite{4JL,4TN}. Korteweg suggests a modification of the
Navier-Stokes system to account additionally for phase transition
phenomena in introducing a term of capillarity. He assumed that the
thickness of the interfaces was not null as in the {\it sharp
interface approach}. This is called the {\it diffuse interface
approach}.\\
Korteweg-type models are based on an extended version of
nonequilibrium thermodynamics, which assumes that the energy of the
fluid not only depends on standard variables but on the gradient of
the density. In terms of the free energy, this principle takes the
form of a generalized Gibbs relation, see \cite{4TN}.\\
In the present chapter, we follow a new approach introduced by
Coquel, Rohde and their collaborators in \cite{4CR}. They remark
that the local diffuse interface approach requires more regular
solution than in the original sharp interface approach. Indeed the
interfaces are assumed of non zero thickness, so that the density
vary continuously between the two interfaces, whereas in the sharp
interface models, the interfaces represent zone of discontinuity for
the density. Coquel, Rohde and their collaborators present an
alternative model with a capillarity term which avoids spatial
derivatives. The model reads:
$$
\begin{aligned}
\begin{cases}
&\p_{t}\rho+{\rm div}(\rho u)=0\\
&\p_{t}(\rho\,u)+{\rm div}(\rho u\otimes u)-\mu\D u-(\lambda+\mu)\n{\rm div} u+\n(P(\rho))=\kappa\rho\n D[\rho]\\
&(\rho_{t=0},u_{t=0})=(\rho_{0},u_{0})\\
\end{cases}
\end{aligned}
\leqno{(NSK)}
$$
with:
$$\mu>0\;\;\mbox{and}\;\;\lambda+2\mu>0$$
where $\rho$ denotes the density of the fluid and $u\in\R^{N}$ the
velocity, $\mu$ and $\lambda$ represent the viscosity coefficients,
$\kappa$ is a coefficient of capillarity, $P$ is a general pressure
function. We are particulary interested by Van der Waals type
pressure:
$$
\begin{aligned}
&P:(0,b)\rightarrow(0,+\infty)\\
&P(\rho)=\frac{RT_{*}\rho}{b-\rho}-a\rho^{2}\\
\end{aligned}
$$
where $a$, $b$, $R$, $T_{*}$ are positive constants, $R$ being the
specific gas constant. For fixed values $a$, $b$ we choose the
constant reference temperature $T_{*}$ so small as $P$ to be
monotone decreasing in
some non-empty interval.\\
Further we impose the conditions:
\begin{equation}
u(t,x)\rightarrow0,\;\;\rho(t,x)\rightarrow
0\;\;\mbox{as}\;\;|x|\rightarrow+\infty, \label{4bord}
\end{equation}
In last section, we consider also more general situations: monotone
pressure, other boundary conditions than boundary condition
\ref{4bord}, at infinity instead of \ref{4bord}:
\begin{equation}
u(t,x)\rightarrow u_{\infty},\;\;\rho(t,x)\rightarrow
\rho_{\infty}\;\;\mbox{as}\;\;|x|\rightarrow+\infty, \label{4bord2}
\end{equation}
where $\rho_{\infty}$ is a given nonnegative constant.\\
The term $\kappa\rho\n D[\rho]$ corresponds to the capillarity which
is supposed to model capillarity effects close to phase transitions
in \cite{4Ko}. The classical Korteweg's capillarity term is
$D[\rho]=\D\rho$.\\
Based  on Korteweg's original ideas Coquel, Rohde and their
collaborators in \cite{4CR} and Rohde in \cite{4Ro} choose a
nonlocal capillarity term $D$ which penalizes rapid variations in
the density field close from the interfaces. They introduce the
following capillarity term:
$$D[\rho]=\phi*\rho-\rho$$
where $\phi$ is chosen so that:
$$\phi\in L^{\infty}(\R^{N})\cap C^{1}(\R^{N})\cap W^{1,\,1}(\R^{N}),\;\;\;\int_{\R^{N}}\phi(x)dx=1,\;\;\phi
\;\;\mbox{even},\;\mbox{and}\;\;\phi\geq0.
$$
This choice of capillarity term allows to get solution with jumps,
i.e with sharp interfaces. Before tackling the global stability
theory for the system $(NSK)$, let us derive formally the uniform
bounds available on $(\rho,u)$.
\subsection{Energy spaces}
We assume to simplify that $P(\rho)=a\rho^{\gamma}$ with
$\gamma\geq1$. Let $\Pi $ (free energy) be defined by:
\begin{equation}
\Pi(s)=s\biggl(\int^{s}_{0}\frac{P(z)}{z^{2}}dz\biggl),
\label{4pressionpi}
\end{equation}
so that $P(s)=s\Pi^{'}(s)-\Pi(s)\, ,\,\Pi^{'}(\bar{\rho})=0$ and if
we renormalize the mass equation:
$$\p_{t}\Pi(\rho)+{\rm div}(u\Pi(\rho))+P(\rho){\rm
div}(u)=0\;\;\mbox{in}\;\;{\cal D}^{'}((0,T)\times\R^{N}).$$
 Notice that
$\Pi$ is convex.
Multiplying the equation of momentum conservation by $u$ and integrating by parts over $\R^{N}$,
we obtain the following energy estimate:\\
\begin{equation}
\begin{aligned}
&\int_{\R^{N}}(\frac{1}{2}\rho
|u|^{2}+\Pi(\rho)+E_{global}[\rho(.,t)])(x)dx(t)+\int_{0}^{t}\int_{\R^{N}}(\mu
D(u):D(u)\\
&\hspace{2cm}+(\lambda+\mu)|{\rm div} u|^{2})dx
\leq\int_{\R^{N}}\big(\frac{|m_{0}|^{2}}{2\rho}+\Pi(\rho_{0})+E_{global}[\rho_{0}]\big)dx,
\end{aligned}
\label{4inegaliteenergie}
\end{equation}
where we have:
$$E_{global}[\rho(.,t)](x)=\frac{\kappa}{4}\int_{\R^{N}}\phi(x-y)(\rho(y,t)-\rho(x,t))^{2}dy.$$
The only non-standard term is the energy term $E_{global}$ which
comes from the product of $u$ with the capillarity term
$\kappa\rho\n(\phi*\rho-\rho)$. Indeed we have:
$$
\begin{aligned}
\kappa\int_{\R^{N}}u(t,x)\rho(t,x)\cdot\n&([\phi*\rho(t,\cdot)](x)-\rho(t,x))dx\\
&\;\;\;=-\kappa\int_{\R^{N}}{\rm div}(u(t,x)\rho(t,x))([\phi*\rho(t,\cdot)](x)-\rho(t,x))dx,\\
&\;\;\;=\kappa\int_{\R^{N}}\frac{\p}{\p t}\rho(t,x)([\phi*\rho(t,\cdot))](x)-\rho(t,x))dx,\\
&\;\;\;=-\frac{d}{dt}\int_{\R^{N}}E_{global}[\rho(t,\cdot)](x)dx\;.\\
\end{aligned}
$$
To derive the last equality we use the relation:
$$
\begin{aligned}
\frac{d}{dt}\int_{\R^{N}}E_{global}[\rho(t,\cdot)](x)dx&=\frac{\kappa}{2}
\int_{\R^{N}}\int_{\R^{N}}\phi(x-y)(\rho(t,y)-\rho(t,x))\frac{\p}{\p
t}\rho(t,y)dydx\\
&\;\;\;\;\;\;\;\;\;\;\;+\frac{\kappa}{2}\int_{\R^{N}}\int_{\R^{N}}\phi(y-x)(\rho(t,x)-
\rho(t,y))\frac{\p}{\p
t}\rho(t,x)dydx,\\[2mm]
&=\kappa\int_{\R^{N}}\int_{\R^{N}}\phi(x-y)(\rho(t,y)-\rho(t,x))\frac{\p}{\p
t}\rho(t,y)dydx,\\[2mm]
&=-\kappa\int_{\R^{N}}([\phi*\rho(t,\cdot)](x)-\rho(t,x))\frac{\p}{\p
t}\rho(t,x)dx\;.\\
\end{aligned}
$$
where we just use integration by parts.\\
In the sequel we will note:
\begin{equation}
{\cal E}(\rho,\rho u)(t)=\int_{\R^{N}}(\frac{1}{2}\rho
|u|^{2}+\Pi(\rho)+E_{global}[\rho(.,t)])(x)dx(t),
\label{4defenergie}
\end{equation}
We are interested to use the above inequality energy to determine the functional space we must work with.\\
So if we expand  $E_{global}[\rho(.,t)](x)$ we get:\\
$$E_{global}[\rho(t,\cdot)](x)=\frac{\kappa}{4}\big(\rho^{2}+\phi*\rho^{2}-2\rho\,(\phi*\rho)\big).$$
Because by the mass equation we obtain that $\rho$ is bounded in $
L^{\infty}(0,T;L^{1}(\R^{N}))$ if we suppose that $\rho_{0}\in
L^{1}$ and we have supposed that $\phi\in L^{\infty}(\R^{N})$, we
obtain that $\rho\,(\phi*\rho)$ is bounded in $
L^{\infty}(0,T;L^{1}(\R^{N}))$. So we get that
$\rho^{2}+\phi*\rho^{2}\in L^{\infty}(0,T;L^{1}(\R^{N}))$ and as
$\phi\geq 0$ and $\rho\geq0$ we get a control of $\rho$ in
$L^{\infty}(0,T;L^{2}(\R^{N}))$ (a property which turns out to be
important to taking advantage of the theory of renormalized
solutions, indeed $\rho$ in $L^{\infty}(0,T;L^{2}(\R^{N}))$ implies
that $\rho\in L^{2}_{loc}(\R^{+}\times\R^{N})$ and we will can use
the theorem of Diperna-Lions on renormalized solutions, see
\cite{4L1})
\\
In view of (\ref{4inegaliteenergie}), we can specify initial
conditions on $\rho_{/t=0}=\rho_{0}$ and $\rho u_{/t=0}=m_{0}$ where
we assume that:
\begin{equation}
\begin{aligned}
&\bullet\;\;\rho_{0}\geq 0\;\;\mbox{a.e in}\;\;\R^{N},\;\rho_{0}\in
L^{1}(\R^{N})\cap L^{s}(\R^{N})\;\;\;\mbox{with}\;\;s=\max(2,\gamma),\hspace{6cm}\\
&\bullet\;\;m_{0}=0\;\;\mbox{a.e on}\;\;{\rho_{0}=0},\\
&\bullet\;\;\frac{|m_{0}|^{2}}{\rho_{0}}\;\;\mbox{(defined to be $0$
on $
{\rho_{0}=0}$) is in}\;\;L^{1}(\R^{N}).\\
\end{aligned}
\label{4donneeinitiale}
\end{equation}
We deduce the following a priori bounds which give us the energy
space in which we will work:
\begin{itemize}
\item $\rho\in L^{\infty}(0,T;L^{1}(\R^{N})\cap L^{s}(\R^{N}))$,
\item $\rho|u|^{2}\in L^{\infty}(0,T;L^{1}(\R^{N}))$,
\item $\n u\in L^{2}((0,T)\times\R^{N})^{N}$.
\end{itemize}
We will use this uniform bound in our result of compactness. Let us
emphasize at this point that the above a priori bounds do not
provide any control on $\n\rho$  in contrast with the case of
$D[\rho]=\D\rho$ (see in \cite{4DD}).\\
\subsection{Notion of weak solutions}
We now explain what we mean by
renormalized weak solutions, weak
solutions, and bounded energy weak solution of problem $(NSK)$.\\
Multiplying mass equation by $b^{'}(\rho)$, we obtained the
so-called renormalized equation (see \cite{4L1}):
\begin{equation}
\frac{\p}{\p t}b(\rho)+{\rm div}(b(\rho)u)+(\rho
b^{'}(\rho)-b(\rho)){\rm div}u=0. \label{4renormalis}
\end{equation}
with:
\begin{equation}
b\in C^{0}([0,+\infty))\cap C^{1}((0,+\infty)),\;\;|b^{'}(t)|\leq
ct^{-\lambda_{0}},\;\;t\in(0,1],\;\,\lambda_{0}<1 \label{4foncb1}
\end{equation}
with growth conditions at infinity:
\begin{equation}
|b^{'}(t)|\leq
ct^{\lambda_{1}},\;\;t\geq1,\;\;\mbox{where}\;\;c>0,\;-1<\lambda_{1}<\frac{s}{2}-1.
\label{4foncb2}
\end{equation}
\begin{definition}
A couple $(\rho,u)$ is called a renormalized weak solution of
problem $(NSK)$ if we have:
\begin{itemize}
\item Equation of mass holds in ${\cal D}^{'}(\R^{N})$.
\item Equation (\ref{4renormalis}) holds in ${\cal D}^{'}(\R^{N})$ for any function $b$ belonging to (\ref{4foncb1})
and (\ref{4foncb2}).
\end{itemize}
\label{4defrenormal}
\end{definition}
\begin{definition}
Let the couple $(\rho_{0},u_{0})$ satisfy;
\begin{itemize}
\item $\rho_{0}\in L^{1}(\R^{N})$, $\Pi(\rho_{0})\in L^{1}(\R^{N})$, $E_{global}[\rho_{0}]\in L^{1}(\R^{N})$,
$\rho_{0}\geq 0$ a.e in $\R^{N}$.
\item $\rho_{0}u_{0}\in (L^{1}(\R^{N}))^{N}$ such that $\rho_{0}|u_{0}|^{2}1_{\rho_{0}>0}\in L^{1}(\R^{N})$
\item and such that $\rho_{0}u_{0}=0$ whenever $x\in\{\rho_{0}=0\}$,
\end{itemize}
where the quantity $\Pi$ is defined in (\ref{4pressionpi}). We have
the following definitions:
\begin{enumerate}
\item A couple $(\rho,u)$ is called a weak solution of problem
$(NSK)$ on $\R$ if:
\begin{enumerate}
\item $\rho\in L^{r}(L^{r}(\R^{N})$ for $s\leq r\leq+\infty$,
\item $P(\rho)\in L^{\infty}(L^{1}(\R^{N}))$, $\rho\geq0$ a.e in
$\R\times\R^{N}$,
\item $\n u\in L^{2}(L^{2}(\R^{N}))$, $\rho|u|^{2}\in
L^{\infty}(L^{1}(\R^{N}))$.
\item Mass equation holds in ${\cal D}^{'}(\R\times\R^{N})$.\label{4vrai1}
\item Momentum equation holds in ${\cal D}^{'}(\R\times\R^{N})^{N}$.\label{4vrai2}
\item $\lim_{t\rightarrow
0^{+}}\int_{\R^{N}}\rho(t)\va=\int_{\R^{N}}\rho_{0}\va$,
$\forall\va\in{\cal D}(\R^{N})$,\label{4vrai3}
\item $\lim_{t\rightarrow
0^{+}}\int_{\R^{N}}\rho u(t)\cdot\phi=\int_{\R^{N}}\rho_{0}
u_{0}\cdot\phi$, $\forall\phi\in{\cal D}(\R^{N})^{N}$.\label{4vrai4}
\end{enumerate}
\item A couple $(\rho,u)$ is called a bounded energy weak solution
of problem $(NSK)$ if in addition to (\ref{4vrai1}), (\ref{4vrai2}),
(\ref{4vrai3}), (\ref{4vrai4}) we have:
\begin{itemize}
\item The quantity ${\cal E}_{0}$ is finite and inequality (\ref{4inegaliteenergie}) with
${\cal E}$ defined by (\ref{4defenergie}) and with ${\cal E}_{0}$ in
place of ${\cal E}(\rho(0),\rho u(0))$ holds a.e in $\R$.
\end{itemize}
\end{enumerate}
\label{4defsolutionsfaibles}
\end{definition}
\subsection{Mathematical results}
We wish to prove global stability results for $(NSK)$ with
$D[\rho]=\phi*\rho-\rho$ in functional spaces very close to energy
spaces. In the non capillary case and $P(\rho)=a\rho^{\gamma}$, P-L.
Lions in \cite{4L2} proved the global existence of weak solutions
$(\rho,u)$ to $(NSK)$ with $\kappa=0$ (which become the system of
Navier-Stokes compressible isotherm) for $\gamma> \frac{N}{2}$ if
$N\geq 4$, $\gamma\geq \frac{3N}{N+2}$ if $N=2,3$ and initial data
$(\rho_{0},m_{0})$ such that:
$$\rho_{0},\;\;\rho_{0}^{\gamma},\;\;\frac{|m_{0}|^{2}}{\rho_{0}}\in
L^{1}(\R^{N}).$$ where we agree that $m_{0}=0$ on
$\{x\in\R^{N}/\;\rho_{0}(x)=0\}$. More precisely, he obtains the
existence of global weak solutions $(\rho,u)$ to $(NSK)$ with
$\kappa=0$ such that for all $t\in(0,+\infty)$:
\begin{itemize}
\item $\rho\in L^{\infty}(0,T;L^{\gamma}(\R^{N}))$ and $\rho\in C([0,T],L^{p}(\R^{N}))$ if $1\leq p<\gamma$,
\item $\rho\in L^{q}((0,T)\times\R^{N})$ for $q=\gamma-1+\frac{2\gamma}{N}>\gamma$.
\item $\rho|u|^{2}\in L^{\infty}(0,T;L^{1}(\R^{N}))$ and $Du\in L^{2}((0,T)\times\R^{N})$.
\end{itemize}
Notice that the main difficulty for proving Lions' theorem consists
in exhibiting strong compactness properties of the density $\rho$ in
$L^{p}_{loc}(\R^{+}\times\R^{N})$ spaces required to pass to the
limit in the pressure
term $P(\rho)=a\rho^{\gamma}$.\\
Let us mention that Feireisl in \cite{4F} generalized the result to
$\gamma>\frac{N}{2}$ in establishing that we can obtain renormalized
solution without imposing that $\rho\in
L^{2}_{loc}(\R^{+}\times\R^{N})$ (what needed Lions in dimension
$N=2,3$, that's why $\gamma-1+\frac{2\gamma}{N}\geq2$), for this he
introduces the concept of oscillation defect measure evaluating the
loss of compactness. We refer to the book of Novotn\'y and
Stra$\check{\mbox{s}}$kraba for more details
(see \cite{4NS}).\\
Let us mention here that the existence of strong solution with
$D[\rho]=\D\rho$ is known since the work by Hattori an Li in
\cite{4Ha1}, \cite{4Ha2} in the whole space $\R^{N}$. In \cite{4DD},
Danchin and Desjardins study the well-posedness of the problem for
the isothermal case with constant coefficients in critical Besov
spaces. We recall too the results from Rohde in \cite{4Ro} who
obtains the existence and uniqueness in finite time for
two-dimensional initial data in $H^{4}(\R^{2})\times H^{4}(\R^{2})$. \\
In the present chapter, we aim at showing the global stability of
weak solutions in the energy spaces for the system $(NSK)$. This
work is composed of four parts, the first one concerns estimates on
the density to get a gain of integrability on the density needed to
pass to the weak limit in the term of pressure and of capillarity.
The second part is the passage to the weak limit in the non-linear
terms of the density and the velocity according to Lions' methods.
The idea is to use renormalized solution to test the weak limit on
convex test functions. In this part we will concentrate on the case
of simple pressure of type $P(\rho)=a\rho^{\gamma}$. We get the
following theorem where
$(\rho_{n},u_{n})_{n\in\mathbb{N}}$ is a suit of global regular
approximate solution for the problem $(NSK)$.
\begin{theorem}
\label{4principal}
Let $N\geq2$. Let $\gamma>N/2$ if $N\geq4$ and $\gamma\geq1$ else.\\
Let the couple $(\rho_{0}^{n},u_{0}^{n})$ satisfy:
\begin{itemize}
\item $\rho^{0}_{n}$ is uniformly bounded in $L^{1}(\R^{N})\cap
L^{s}(\R^{N})$ with $s=\max (\gamma,2)$ and $\rho_{0}^{n}\geq 0$ a.e
in $\R^{N}$,
\item $\frac{|\rho^{0}_{n}u_{0}^{n}|^{2}}{\rho^{0}_{n}}$ is uniformly bounded in
$L^{1}(\R^{N})$,
\item and such that $\rho_{0}^{n}u_{0}^{n}=0$ whenever $x\in\{\rho_{0}=0\}$.
\end{itemize}
In addition we suppose that $\rho^{0}_{n}$ converges in
$L^{1}(\R^{N})$ to $\rho_{0}$, then up to a subsequence
$(\rho_{n},u_{n})$ converges strongly to a weak solution $(\rho,u)$
of the system $(NSK)$ satisfying the initial condition
$(\rho_{0},u_{0})$ as in (\ref{4donneeinitiale}). Moreover we have
the following convergence:
\begin{itemize}
\item $\rho_{n}\rightarrow_{n}\rho$ in $C([0,T],L^{p}(\R^{N}))\cap
L^{r}((0,T)\times K)$ for all $1\leq p<s$, $1\leq r<q$, with
$q=s+\frac{N\gamma}{2}-1$ if $N\geq3$, where $K=\R^{N}$ except when
$N\geq 4$ and $\gamma<\frac{N}{2}(1+\frac{1}{N})$ where $K$ is an
arbitrary compact.
\item $\rho_{n}\rightarrow_{n}\rho$ in $C([0,T],L^{p}(\R^{2}))\cap
L^{r}((0,T)\times K)$ for all $1\leq p<s$, $1\leq r<q$, with $K$ an
arbitrary compact in $\R^{2}$ if $N=2$.
\end{itemize}
In addition we have:
\begin{itemize}
\item $\rho_{n}u_{n}\rightarrow\rho
u$ in $L^{p}(0,T;L^{r}(\R^{N}))$ for all $1\leq p<+\infty$ and
$1\leq r<\frac{2s}{s+1}$,
\item
$\rho_{n}(u_{i})_{n}(u_{j})_{n}\rightarrow\rho_{n}u_{i}u_{j}$ in
$L^{p}(0,T;L^{1}(\Omega))$ for all $1\leq p<+\infty$, $1\leq i,j\leq
N$ if $N\geq3$,  where $K=\R^{N}$ except when $N\geq 4$ and
$\gamma<\frac{N}{2}(1+\frac{1}{N})$ where $K$ is an arbitrary
compact.
\item
$\rho_{n}(u_{i})_{n}(u_{j})_{n}\rightarrow\rho_{n}u_{i}u_{j}$ in
$L^{p}(0,T;L^{1}(\Omega))$ for all $1\leq p<+\infty$, $1\leq i,j\leq
N$ with $\Omega$ an arbitrary bounded open set in $\R^{2}$ if $N=2$.
\end{itemize}
\end{theorem}
In the third part we will focus on general pressure, and particulary
van der Waal's pressure. In the fourth part we concentrate on the
case with initial data close to a constant $\bar{\rho}$, and we will
work in Orlicz space, this case is the most adapted for the strong
solution because it enables us to control the vacuum so that one can
use the property of ellipticity of the momentum equation.
\section{Existence of weak solution for a isentropic pressure law}
\subsection{A priori estimates on the density}
In this part we are interested by getting a gain of integrability on
the density and we consider the case where $P(\rho)=a\rho^{\gamma}$.
This will enable us to pass to the weak limit in the pressure and
the Korteweg terms. It is expressed by the following theorem:
\begin{theorem}
\label{4integration} Let $N\geq2$ and $\gamma\geq1$, with in
addition $\gamma>\frac{N}{2}$ if $N\geq4$. Let $(\rho,u)$ be a
regular solution of the system $(NSK)$ with $\rho\geq 0$ and
$\rho\in L^{\infty}(L^{1}\cap L^{s+\e})$ where we define $\e$ below.
Then we have if $\gamma\geq \frac{N}{2}(1+\frac{1}{N})$ for $N\geq4$
:
$$
\begin{aligned}
\int_{(0,T)\times\R^{N}}\big(\rho^{\gamma+\e}+\rho^{2+\e}\big)dxdt\leq
M&\;\;\mbox{for any}\;\;0<\e\leq \frac{2}{N}\gamma-1\;\;\;\mbox{if}
\;\;N\geq4,\\
&\;\;\mbox{and}\;\;\;0<\e\leq\frac{4}{N}-1\;\;\;\mbox{if}\;\; N=2,3.
\end{aligned}
$$
with $M$ depending only on the initial conditions and on the time
$T$.\\
If $\gamma<\frac{N}{2}(1+\frac{1}{N})$ for $N\geq4$, we have:
$$
\begin{aligned}
\int_{(0,T)\times K}\big(\rho^{\gamma+\e}+\rho^{2+\e}\big)dxdt\leq
M^{'}&\;\;\mbox{for any}\;\;0<\e\leq
\frac{2}{N}\gamma-1\;\;\;\mbox{if}
\;\;N\geq4,\\
\end{aligned}
$$
for any arbitrary compact $K$ with $M^{'}$ depending only on the
initial conditions, on $K$ and on the time $T$.
\end{theorem}
{\bf Proof:}
\\
\\
We will begin with the case where $N\geq3$ and we treat after the
specific case $N=2$.
\subsubsection*{Case $N\geq3$:}
We apply to the momentum equation the operator $(-\D)^{-1}{\rm div}$
in order to concentrate us on the pressure and we get:
\begin{equation}
\begin{aligned}
a\rho^{\gamma}=\frac{\p}{\p t}(-\D)^{-1}{\rm div}(\rho u)+
(-\D)^{-1}\p^{2}_{i,j}(\rho u_{i}u_{j})+&(2\mu+\lambda){\rm div} u\\
&-\kappa(-\D)^{-1}{\rm div}\big(\rho\n(\phi*\rho-\rho)\big),\\
\end{aligned}
\label{4premiereeq}
\end{equation}
and in multiplying by $\rho^{\e}$ with $0<\e\leq\min
(\frac{1}{N},\frac{2}{N}\gamma-1)$ to estimate $\rho^{\gamma+\e}$,
we get:
\begin{equation}
\begin{aligned}
a\rho^{\gamma+\e}+\frac{\kappa}{2}\rho^{2+\e}=&-\kappa\rho^{\e}(-\D)^{-1}{\rm
div}(\rho(\n\phi*\rho))+\rho^{\e}(-\D)^{-1}\p^{2}_{ij}(\rho
u_{i}u_{j})\\
&+\frac{\p}{\p t}\big(\rho^{\e}(-\D)^{-1}{\rm div}(\rho
u)\big)-[\frac{\p}{\p t}\rho^{\e}](-\D)^{-1}{\rm div}(\rho
u)+(\mu+\zeta){\rm div} u\;,\\
\end{aligned}
\end{equation}
where we note $\xi=\lambda+\mu$. We now rewrite the previous
equality as follows:
\begin{equation}
\begin{aligned}
&a\rho^{\gamma+\e}+\frac{\kappa}{2}\rho^{2+\e}=-\kappa\rho^{\e}(-\D)^{-1}{\rm
div}(\rho(\n\phi*\rho))+\rho^{\e}(-\D)^{-1}\p^{2}_{ij}(\rho
(u_{i})(u_{j}))\\
&\hspace{1,3cm}+\frac{\p}{\p t}\big(\rho^{\e}(-\D)^{-1}{\rm
div}(\rho u)\big)+{\rm div}[u\rho^{\e}(-\D)^{-1}{\rm div}(\rho
u)]+(\mu+\zeta){\rm div}
u\\
&\hspace{2,8cm}-(\rho)^{\e}u\cdot\n(-\D)^{-1}{\rm div}(\rho
u)+(1-\e)({\rm
div}u)\rho^{\e}(-\D)^{-1}{\rm div}(\rho u)\;,\\
\end{aligned}
\label{42.14}
\end{equation}
\\
Next we integrate (\ref{42.14}) in time on $[0,T]$ and in space so
we obtain:
\begin{equation}
\begin{aligned}
&\int_{(0,T)\times\R^{N}}(a\rho^{\gamma+\e}+\frac{\kappa}{2}\rho^{2+\e})dx\,dt=
\int_{(0,T)\times\R^{N}}\biggl(\frac{\p}{\p
t}[\rho^{\e}(-\D)^{-1}{\rm div}(\rho u)]+
(\mu+\zeta)({\rm div}u)\rho^{\e}\\[2mm]
&\hspace{2,9cm}+(1-\e)({\rm div} u)\rho^{\e}(-\D)^{-1}{\rm div}(\rho
u)+\rho^{\e}[R_{i}R_{j} (\rho u_{i}u_{j})-u_{i}R_{i}R_{j}(\rho
u_{j})]\\[2mm]
&\hspace{3,8cm}+{\rm div}[u\rho^{\e}(-\D)^{-1}{\rm div}(\rho
u)]-\kappa\rho^{\e}(-\D)^{-1}
{\rm div}(\rho\n(\phi*\rho))\biggl)dx\,dt\;.\\[2mm]
\end{aligned}
\label{4eq3}
\end{equation}
where $R_{i}$ is the classical Riesz transform.\\
Now we want to control the term
$\int^{T}_{0}\int_{\R^{N}}\big(\rho^{\gamma+\e}+\frac{\kappa}{2}\rho^{2+\e}\big)dxdt$.
As $\rho$ is positive, it will enable us to control
$\|\rho\|_{L_{t,x}^{\gamma+\e}}$ and $\|\rho\|_{L_{t,x}^{2+\e}}$.
This may be achieved by  controlling each term on the right
side of (\ref{4eq3}).\\
\\
We start by treating the term $\int_{(0,T)\times\R^{N}}\frac{\p}{\p
t}[\rho^{\e}(-\D)^{-1}{\rm div}(\rho u)]$. So we need to control
$\rho^{\e}(-\D)^{-1}{\rm div}(\rho u)$ in
$L^{\infty}(0,T;L^{1}(\R^{N}))$ and $\rho^{\e}_{0}(-\D)^{-1}{\rm
div}(\rho_{0}u_{0})$ because:
$$
\begin{aligned}
\int_{(0,T)\times\R^{N}}\frac{\p}{\p t}[\rho^{\e}(-\D)^{-1}{\rm
div}(\rho u)](t,x)dt\,dx=&\int_{\R^{N}}\rho^{\e}(-\D)^{-1}{\rm
div}(\rho u)](x)dx(t)\\
&\hspace{1cm}-\int_{\R^{N}}\rho_{0}^{\e}(-\D)^{-1}{\rm div}(\rho_{0}
u_{0})](x)dx,\\
\end{aligned}
$$ We recall that $\rho$,
$\rho^{2}$, $\rho^{\gamma}$ and $\rho |u|^{2}$ are bounded in
$L^{\infty}(L^{1})$ while $Du$ is bounded in
$L^{2}((0,T)\times\R^{N})$ and $u$ is bounded in
$L^{2}(0,T;L^{\frac{2N}{N-2}}(\R^{N}))$  by Sobolev embedding. In
particular by H\"older inequalities we get that $\rho u$ is bounded
in $L^{\infty}(0,T,(L^{\frac{2\gamma}{\gamma+1}}\cap
L^{\frac{4}{3}})(\R^{N}))$. Thus we get in using H\"older inequalities and Sobolev embedding:\\
$\rho^{\e}(-\D)^{-1}{\rm div}(\rho u)\in L^{\infty}(0,T,L^{1} \cap
L^{\alpha})$ with:
$$\frac{1}{\alpha}=\frac{\e}{s}+\min(\frac{\gamma+1}{2\gamma},\frac{3}{4})-\frac{1}{N}
<1.$$
The fact that
$\rho^{\e}(-\D)^{-1}{\rm div}(\rho u)\in L^{\infty}(0,T,L^{1})$ is
obtained by interpolation because $\rho\in L^{\infty}(L^{1})$ and in
using less integrability in Sobolev embedding.
\\
Next we have the same type of estimates for
$\|\rho^{\e}_{0}(-\D)^{-1}{\rm
div}(\rho_{0}u_{0})\|_{L^{1}(\R^{N})}$.\\
Finally (\ref{4eq3}) is rewritten on the following form in using
Green formula:
\\
$$
\begin{aligned}
\int^{T}_{0}\int_{\R^{N}}\big(\rho^{\gamma+\e}&+\frac{\kappa}{2}\rho^{2+\e}\big)dxdt\leq
C\big(1+
\int^{T}_{0}\int_{\R^{N}}\big[\,|{\rm div}u|\rho^{\e}(1+|(-\D)^{-1}{\rm div}(\rho u)|\,]\\[3mm]
&+\rho^{\e}|R_{i}R_{j}(\rho u_{i}u_{j})-u_{i}R_{i}R_{j}(\rho u_{j})|
+\kappa\rho^{\e}|(-\D)^{-1}{\rm div}(\rho\n(\phi*\rho))\big]dt\,dx\big).\\
\end{aligned}
$$
Now we will treat each term of the right hand side. We treat all the
terms with the same type of estimates than P.-L. Lions in \cite{4L2}
, excepted the capillarity term.
\\
\\
We start with the term $|{\rm div}u|\rho^{\e}|(-\D)^{-1}{\rm
div}(\rho u)|$ where we have:
$$|{\rm div}u|\in L^{2}(L^{2}),\;\;
\rho^{\e}\in L^{\infty}(L^{\frac{s}{\e}}),\;\; \rho\,u\in
L^{2}(0,T,L^{r}(\R^{N}))$$ with
$\frac{1}{r}=\frac{1}{s}+\frac{N-2}{2N}$ and by Sobolev embedding
$|(-\D)^{-1}{\rm div}(\rho u)|\in L^{2}(L^{s^{'}})$ with
$\frac{1}{s^{'}}=\frac{1}{r}-\frac{1}{N}$ (this is possible only if
$r<N$). We are in a critical case for the Sobolev embedding ( i.e
$r\geq N$) only when $N=3$ and $\gamma\geq 6$, that's
why for $N=3$ and $\gamma\geq6$.\\
So by H\"older inequalities we get $|{\rm
div}u|\rho^{\e}|(-\D)^{-1}{\rm div}(\rho u)|\in L^{1}(L^{s_{1}})$
with:
$\frac{1}{s_{1}}=\frac{1}{s^{'}}+\frac{\e}{s}+\frac{1}{2}=1-\frac{2}{N}+\frac{1+\e}{s}\leq1$
as we have $s>\frac{N}{2}$.\\
Moreover by interpolation $|{\rm div}u|\rho^{\e}|(-\D)^{-1}{\rm
div}(\rho u)|$ belongs to $L^{1}(0,T;L^{1}(\R^{N}))$.\\
We now treat the case $N=3$ and $\gamma\geq6$ where we choose in
this case $\e=\frac{2}{N}\gamma-1$ to explicit precisely this case:
$$
\begin{aligned}
&\||{\rm div}u|\rho^{\e}|(-\D)^{-1}{\rm div}(\rho u)|\|_{L^{1}}\leq
\|Du\|_{L^{2}(L^{2})}\|\rho\|_{L^{\gamma+\e}}^{\e}\|\rho
u\|_{L^{\frac{2(\gamma+\e)}{\gamma-\e}}
(L^{\frac{6(\gamma+\e)}{5\gamma-\e}})}\\
&\leq C\|\rho\|^{\e}_{L^{\gamma+\e}}\|\rho
u\|_{L^{\frac{10\gamma-6}{\gamma+3}}
(L^{\frac{3(10\gamma-6)}{13\gamma+3}})}\leq\, C\|\rho\|^{\e}_
{L^{\gamma+\e}}\|\rho u\|^{\frac{\gamma+3}{5\gamma-3}}_{L^{2}
(L^{\frac{6\gamma}{\gamma+6}})}\|\rho u\|^{\frac{2(2\gamma-3}{5\gamma-3}}_{L^{\infty}(L^{2})}\\
&\leq C\|\rho\|^{\e}_{L^{\gamma+\e}}\|\rho
u\|_{L^{2}(L^{\frac{6\gamma}{\gamma+6}})}^{\frac{5\gamma}
{5\gamma-3}}\|\rho u\|_{L^{\infty}(L^{\frac{2\gamma}{\gamma+1}})}^{\frac{2(2\gamma-3)}{5\gamma-3)}}\\
&\leq C\|\rho\|^{\e}_{L^{\gamma+\e}}
\end{aligned}
$$
\\
since we have
$\frac{1}{2}+\frac{\e}{\gamma+\e}+\frac{\gamma-\e}{2(\gamma+\e)}=1$,
$\frac{1}{2}+\frac{\e}
{\gamma+\e}+\frac{5\gamma-\e}{6(\gamma+\e)}-\frac{1}{3}=1$, and
$\frac{6(\gamma+\e)}{5\gamma-\e}=3\frac{10\gamma-6}
{13\gamma+3}<3$.\\
\\
We now want to treat the term: $\rho^{\e}|(-\D)^{-1}{\rm div
}(\rho\n(\phi*\rho))|$, so we have:
$\rho\n(\phi*\rho)=\rho(\n\phi*\rho)\in L^{\infty}(L^{1}\cap
L^{\frac{s}{2}})$ by H\"older inequalities
and the fact that we have $\rho\in L^{\infty}(L^{1})$ and $\n\phi\in L^{1}$.\\
After we get that $\rho^{\e}(-\D)^{-1}{\rm
div}(\rho\n(\phi*\rho))\in L^{\infty}(L^{r_{1}})$ with:
$\frac{1}{r_{1}}=\frac{\e}{s}+\frac{2}{s}-\frac{1}{N}=\frac{2+\e}{s}-\frac{1}{N}<1$.\\
We conclude that $\rho^{\e}(-\D)^{-1}\rm div(\rho\n(\phi*\rho))$ is
$L^{\infty}(L^{1})$ in using interpolation when $N=2,3$. Indeed we
have $\rho\n(\phi*\rho)\in L^{\infty}(L^{1})$ and in choosing
$\e=\frac{2}{N}s-1$ we have: $1-\frac{1}{N}+\frac{2}{N}s-1\geq1$
when $s\geq\frac{N}{2}(1+\frac{1}{N})$. This is the case
when $N=2,3$, and this is the case when $N\geq 4$ and $\gamma\geq\frac{N}{2}(1+\frac{1}{N})$.\\
In the other case we need to work in arbitrary compact.
\\
We have after the term $({\rm div}(u))\rho^{\e}$. We recall that
$\rho^{\e}$ is in $L^{\infty}(L^{\frac{1}{\e}}\cap
L^{\frac{s}{\e}})$. If $\e\geq\frac{1}{2}$ (i.e
$s\geq\frac{3}{4}N$), the bound is obvious because
$\frac{1}{2}+\e\geq1$ and $\frac{1}{2}+\frac{\e}{s}<1$, we can then
conclude by interpolation. On the other hand, this rather simple
term presents a technical difficulty when $\e\leq\frac{1}{2}$ since
we do not know in that case if ${\rm div u}\,\rho^{\e}\in
L^{1}(\R^{N}\times(0,T))$. One way to get round the difficulty is to
multiply (\ref{4premiereeq}) by $\rho^{\e}1_{\{\rho\geq1\}}$. Then
we obtain an estimate on $\rho^{s+\e}1_{(\rho\geq1)}$ in
$L^{1}((0,T)\times\R^{N})$ as $\rho^{\e}1_{\{\rho\geq1\}}|{\rm
div}u|\leq \rho|{\rm div}u|\in L^{1}((0,T)\times\R^{N})$ (where
$\e\leq\frac{1}{2}$) and we can conclude since
$0\leq\rho^{s+\e}1_{\{\rho<1\}}\leq\rho$ on $(0,T)\times\R^{N}$ and
$\rho\in L^{\infty}(L^{1})$.\\
\\
We end with the following term
$\rho^{\e}(R_{i}R_{j}(\rho\,u_{i}u_{j})-u_{i}R_{i}R_{j}(\rho\,u_{j}))$.
In the same way than in the previous inequalities we have
$\rho^{\e}R_{i}R_{j}(\rho\,u_{i}u_{j})$ is bounded in
$L^{1}(0,T;L^{1}(\R^{N}))$. Indeed we have by H\"older inequalities
and the fact that $R_{i}$ is continuous from $L^{p}$ in $L^{p}$ with
$1<p<+\infty$:
$\frac{1}{s}+2\frac{N-2}{2N}+\frac{\e}{s}=1-\frac{2}{N}+\frac{1+\e}{s}\leq1$
(because $s>\frac{N}{2}$). And we conclude by interpolation. We
treat the term $\rho^{\e}u_{i}R_{i}R_{j}(\rho\,u_{j})$ similarly.
\\
\\
We have to treat now the case $N=2$ where we have to modify the
estimates when we are in critical cases for Sobolev embedding.
\subsubsection*{Case $N=2$:}
In the case $N=2$ most of the proof given above stay exact except
for the slightly more delicate terms $\rho^{e}{\rm
div}u|(-\D)^{-1}{\rm div}(\rho u)|$ and
$\rho^{\e}(R_{i}R_{j}(\rho\,u_{i}u_{j})-u_{i}R_{i}R_{j}(\rho\,u_{j}))$.\\
We start with the term $|\rho^{e}{\rm div}u(-\D)^{-1}{\rm div}(\rho
u)|$. In our previous estimate it was possible to use Sobolev
embedding on the term $(-\D)^{-1}{\rm div}(\rho u)|$ only if $r\geq
N$ (see above the notation), so in the case where $N=2$ we are in a
critical case for
the Sobolev embedding when $\gamma\geq2$.\\
This may be overcome by using that, by virtue of Sobolev embedding,
we have:
$$
\begin{aligned}
&\||{\rm div}u|\rho^{\e}|(-\D)^{-1}{\rm div}(\rho u)|\|_{L^{1}}\leq
C\|\rho\|_{L^{\gamma+\e}(L^{\gamma+\e})}^{\e}
\|\rho u\|_{L^{2(\gamma+\e)}(L^{\frac{2(\gamma+\e)}{\gamma+\e+1}})}\\
\end{aligned}
$$
Indeed:
$\frac{1}{2}+\frac{\e}{\gamma+\e}+\frac{\gamma+\e+1}{2(\gamma+\e)}-\frac{1}{2}=\frac{1}{2}+\frac{2\e+1}{2\e+2\gamma}\leq1=
\frac{1}{2}+\frac{\e}{\gamma+\e}+\frac{1}{2(\gamma+\e})\leq1$ and
$\frac{1}{2}+\frac{\e}{\gamma+\e}+\frac{1}{2(\gamma+\e)}\leq1$.
Moreover we have as $\rho u=\sqrt{\rho}\sqrt{\rho}u$
$$\|\rho u\|_{L^{2(\gamma+\e)}(L^{\frac{2(\gamma+\e)}{\gamma+\e+1}})}\leq C\|
\rho\|_{L^{\gamma+\e}(L^{\gamma+\e})}^{\frac{1}{2}}$$
then:
$$\||{\rm div}u|\rho^{\e}|(-\D)^{-1}{\rm div}(\rho u)|\|_{L^{1}(L^{1})}\leq C\|\rho\|^
{\e+\frac{1}{2}}_{L^{\gamma+\e}(L^{\gamma+\e})}.$$
\\
Next we are interested by the term
$\rho^{\e}(R_{i}R_{j}(\rho\,u_{i}u_{j})-u_{i}R_{i}R_{j}(\rho\,u_{j}))$.
We use the fact that $u$ is bounded in $L^{2}(0,T;\dot{H}^{1})$ and
thus in $L^{2}(0,T;BMO)$. Then by the Coifman-Rochberg-Weiss
commutator theorem in \cite{4CRW},
we have for almost all $t\in[0,T]$:\\
$$\|R_{i}R_{j}(\rho\,u_{i}u_{j})-u_{i}R_{i}R_{j}(\rho\,u_{j})\|_{L^{\frac{2(\gamma+\e)}{\gamma+\e+1}}(L^{\frac{2(\gamma+\e)}{\gamma+\e+1}})}
\leq C\|u\|_{L^{2}(
BMO)}\|\rho\,u\|_{L^{2(\gamma+\e)}(L^{\frac{2(\gamma+\e)}{\gamma+\e+1}})}$$
So we have:\\
$$\|\rho^{\e}(R_{i}R_{j}(\rho\,u_{i}u_{j})-u_{i}R_{i}R_{j}(\rho\,u_{j}))\|_{L^{1}}\leq C\|\rho\|^{\e+\frac{1}{2}}
_{L^{\gamma+\e}(L^{\gamma+\e})}$$
\\
In view of the previous inequalities we get finally:
$$\|\rho\|_{L^{\gamma+\e}(L^{\gamma+\e})}^{\gamma+\e}\leq C(1+\|\rho\|_{L^{\gamma+\e}(L^{\gamma+\e})}
^{\frac{1}{2}+\e})$$ and the $L^{\gamma+\e}(L^{\gamma+\e})$ bound on
$\rho$ is proven since $\frac{1}{2}+\e<\gamma+\e$. \hfill{$\Box$}
\\
\subsection{Compactness results for compressible Navier-Stokes
equations of Korteweg type in the case of isentropic pressure} In
the sequel we are not going to treat in details the case with
$N\geq4$ and $\gamma<\frac{N}{2}(1+\frac{1}{N})$, we just remark
that the proof is the same as in the case $N=2$, it suffices to
localize because we can only apply the theorem \ref{4integration} on
the gain of integrability
on any compact $K$.\\
So let follow the theorem \ref{4integration} and assume that
$\gamma>\frac{N}{2}$ if $N\geq4$ and $\gamma\geq1$ such that if
$(\rho,u)$ is a regular solution then $\rho\in
L^{q}((0,T)\times\R^{N})$ with $q=\gamma+1-\frac{2\gamma}{N}$. We
can observe that in this case $q>s=\max(\gamma,2)$. We will see that
it will be very useful in the sequel to justify the passage to the
weak limit in some terms to get a gain of integrability on the
density. Indeed the key point to proving the existence of weak
solutions is the passage to the limit in the term of pressure and in
the term of capillarity
 $\rho\n(\phi*\rho-\rho)$.\\
First, we assume that a sequence $(\rho_{n},u_{n})_{n\in\mathbb{N}}$
of approximate weak solutions has been constructed by a mollifying
process, which have suitable regularity to justify the formal
estimates like the energy estimate and the theorem
\ref{4integration}. $(\rho_{n},u_{n})_{n\in\mathbb{N}}$
has the initial data of the theorem \ref{4principal} with uniform bounds, i.e:
\\
Moreover $(\rho_{n},u_{n})_{n\in\mathbb{N}}$ verifies the energy
inequalities (\ref{4inegaliteenergie}) and the previous theorem
\ref{4integration}, we have then:
\begin{itemize}
\item $\rho_{n}$ is bounded uniformly in $L^{\infty}(0,T;L^{1}\cap
L^{s}(\R^{N}))\cap C([0,T];L^{p}(\R^{N}))$ for $1\leq
p<\max(2,\gamma)$,
\item $\rho_{n}\geq0$ a.e. $\rho_{n}$ is bounded uniformly
in $L^{q}(0,T,\R^{N})$ with $q>s$,
\item $\n u_{n}$ is bounded in $L^{2}(0,T;L^{2}(\R^{N}))$,
$\rho_{n} |u_{n}|^{2}$ is bounded in
$L^{\infty}(0,T;L^{1}(\R^{N}))$,
\item $u_{n}$ is bounded in $L^{2}(0,T;L^{\frac{2N}{N-2}}(\R^{N}))$ for $N\geq3$.
\end{itemize}
And we have in passing to the weak limit in the previous bound in
extracting subsequence if necessary:
\begin{itemize}
\item $\rho_{n}\rightarrow\rho$ weakly in $L^{s}((0,T)\times\R^{N})$,
\item $u_{n}\rightarrow u$ weakly
in $L^{2}(0,T,\dot{H}^{1}(\R^{N}))$,
\item $\rho_{n}^{\gamma}\rightarrow
\overline{\rho^{\gamma}}$ weakly in $L^{r}((0,T)\times\R^{N})$ for
$r=\frac{q}{\gamma}>1$,
\item $\rho_{n}^{2}\rightarrow
\overline{\rho^{2}}$ weakly in $L^{r_{1}}((0,T)\times\R^{N})$ for
$r_{1}=\frac{q}{2}>1$.
\end{itemize}
\begin{notation}
We will always write in the sequel $\overline{B(\rho)}$ to mean the
weak limit of the sequence $B(\rho_{n})$  bounded in appropriate
space that we will precise.
\end{notation}
We recall that the main difficulty will be to pass to the limit in
the pressure term and the capillary term. The idea of the proof will
be to test the convergence of the sequence
$(\rho_{n})_{n\in\mathbb{N}}$ on convex functions $B$ in order to
use their properties of lower semi-continuity with respect to the
weak topology in $L^{1}(\R^{N})$. In this goal we will use the
theory of renormalized solutions introduced by Diperna and Lions in
\cite{DLion}. So we will obtain strong convergence of $\rho_{n}$ in
appropriate spaces.
\subsection{Idea of the proof\label{4idee}}
We here give a sketchy proof of the theorem \ref{4principal}. In
this spirit we can rewrite mass conservation of the regular solution
$(\rho_{n},u_{n})_{n\in\mathbb{N}}$ on the form:
$$\frac{p}{\p t}(B(\rho_{n}))+{\rm div}(u_{n}B(\rho_{n}))=(B(\rho_{n})-\rho_{n} B^{'}(\rho_{n})){\rm
div}u_{n}.$$ Supposing that $B(\rho_{n})$ is bounded in appropriate
space we can pass to the weak limit where we have in the energy
space $\rho_{n}\rightharpoonup\rho$ and $u_{n}\rightharpoonup u$, so
we get:
\begin{equation}
\frac{p}{\p t}(\overline{B(\rho)})+{\rm div}(u
\overline{B(\rho)})=\overline{(B(\rho)-\rho B^{'}(\rho)){\rm div}u}
\label{41}
\end{equation}
We will set: $b(\rho)=B(\rho)-\rho B^{'}(\rho).$\\
Next in seeing the mass equation for approximate solutions
$(\rho_{n},u_{n})_{n\in\mathbb{N}}$ , and in passing directly to the
limit via the weak convergence arguing like P-L. Lions in \cite{4L2}
p 13 we get:
\begin{equation}
\frac{d}{dt}\rho+{\rm div}(\rho u)=0. \label{4a1}
\end{equation}
After we will just have  to verify the passage to the limit for the
product $\rho u$. Next we will use the theorem on the renormalized
solutions of Diperna-Lions in \cite{4L1} on (\ref{4a1}) in recalling
that $\rho\in L^{\infty}(L^{2})$. So we get:
\begin{equation}
\frac{d}{dt}(B(\rho))+{\rm div}(uB(\rho))=b(\rho){\rm div}(u).
\label{42}
\end{equation}
Next we subtract (\ref{41}) to (\ref{42}), so we obtain:
\begin{equation}
\frac{d}{dt}(\overline{B(\rho)}-B(\rho))+{\rm
div}(u(\overline{B(\rho)}-B(\rho)))=\overline{b(\rho){\rm
div}u}-b(\rho){\rm div}u. \label{43}
\end{equation}
Consequently, in order to estimate the difference
$\overline{B(\rho)}-B(\rho)$ which tests the convergence of
$\rho_{n}$, we need to estimate the difference
$\overline{b(\rho){\rm div}(u)}-b(\rho){\rm div}(u)$. We choose then
$B$ a concave function and we have:
$$\overline{B(\rho)}-B(\rho)\leq 0.$$
The goal will be now to prove the reverse inequality in order to
justify that $B(\rho_{n})$ tends to $B(\rho)$ a.e.\\
So now we aim at estimating the difference $\overline{b(\rho){\rm
div}u}-b(\rho){\rm div}u$. This may be achieved by introducing the
effective viscous pressure $\mbox{P}_{eff}=P-(2\mu+\lambda){\rm
div}u $ after D. Hoff in \cite{4H1}, which satisfies some important
properties of weak
convergence.\\
In fact owing to the capillarity term we adapt Hoff's concept to our
equation in setting:
$$\widetilde{P}_{eff}=P+\frac{\kappa}{2}\rho^{2}-(2\mu+\lambda){\rm div}u.$$
{\bf Proof of theorem \ref{4principal}}
\\
\\
We begin with the case $N\geq3$, and next we will complete the proof
by the case $N=2$ in specifying the changes to bring.\\
Before getting into the heart of the proof, we first recall that we
obtain easily the convergence in distribution sense of
$\rho_{n}u_{n}$ to $\rho u$ and $\rho_{n}(u_{n})_{i}(u_{n})_{j}$ to
$\rho u_{i}u_{j}$. We refer to the classical result by Lions (see
\cite{4L2}) or the book of Novotn\'y and Stra$\check{\mbox{s}}$kraba
\cite{4NS}.
\subsubsection*{Case $N\geq3$}
We have seen in the idea of the proof \ref{4idee} that our goal is
to compare $\overline{B(\rho)}$ and $B(\rho)$ for certain concave
functions $B$ . From the mass equation we have obtained:
\begin{equation}
\p_{t}(\overline{B(\rho)}-B(\rho))+{\rm
div}(u(\overline{B(\rho)}-B(\rho)))=\overline{b(\rho){\rm
div}(u)}-b(\rho){\rm div}(u). \label{4a3}
\end{equation}
\\
So before comparing  $\overline{B(\rho)}$ and $B(\rho)$, we have to
investigate the expression $\overline{b(\rho){\rm
div}(u)}-b(\rho){\rm div}(u)$. By virtue of theorem
\ref{4integration} which gives a gain of integrability we can take
the function $B(x)=x^{\e}$ , as we control for $\e$ small enough
$\rho^{s+\e}$.
Our goal now is to exhibit the effective pressure
$\widetilde{P}_{eff}$, and to multiply it by $\rho^{\e}$ to extract
$\overline{{\rm div}u\, b(\rho)}$ . We will see in the sequel how to
compare it with $b(\rho){\rm div}(u)$. So we focus on the
convergence of the pressure and capillarity terms.
\subsubsection*{Control of the term $\overline{{\rm div}u\, b(\rho)}$ }
So we take the ${\rm div}$ on the momentum equation satisfied by the
regular solution. We get:
\begin{equation}
\begin{aligned}
\frac{\p}{\p t}{\rm div}(\rho_{n} u_{n})+\p^{2}_{ij}(\rho_{n}
u^{i}_{n}u^{j}_{n})-\zeta\D{\rm div}
u_{n}+\D(a\rho^{\gamma}_{n})=&\kappa{\rm
div}(\rho_{n}(\n\phi*\rho_{n}))-\frac{\kappa}{2}\D(\rho^{2}_{n}),\\
\label{410.1}
\end{aligned}
\end{equation}
with $\zeta=\lambda+2\mu$. Applying the operator $(-\D)^{-1}$ to
(\ref{410.1}), we obtain:
\begin{equation}
\begin{aligned}
&\frac{\p}{\p t}(-\D)^{-1}{\rm div}(\rho_{n}
u_{n})+(-\D)^{-1}\p^{2}_{ij}(\rho_{n} u^{i}_{n}u^{j}_{n})+[\zeta{\rm
div}
u_{n}-a\rho^{\gamma}_{n}-\frac{\kappa}{2}\rho^{2}_{n}]\\
&\hspace{8cm}=\kappa(-\D)^{-1}{\rm
div}(\rho_{n}(\n\phi*\rho_{n}))\;.\\
\end{aligned}
\label{43}
\end{equation}
After we multiply (\ref{43}) by $\rho_{n}^{\e}$ with $\e$ that we
choose enough small with $\e\in(0,1)$:
\begin{equation}
\begin{aligned}
&[(\mu+\zeta){\rm div}
u_{n}-a\rho_{n}^{\gamma}-\frac{\kappa}{2}\rho_{n}^{2}]\rho_{n}^{\e}=\kappa\rho_{n}^{\e}(-\D)^{-1}{\rm
div}(\rho_{n}(\n\phi*\rho_{n}))\\
&\hspace{4cm}-\rho_{n}^{\e}\frac{\p}{\p t}(-\D)^{-1}{\rm
div}(\rho_{n} u_{n})-\rho_{n}^{\e}(-\D)^{-1}\p_{ij}(\rho_{n}
u^{i}_{n}u^{j}_{n}).\\
\end{aligned}
\label{44}
\end{equation}
So if we rewrite (\ref{44}), we have:
\begin{equation}
\begin{aligned}
&[(\mu+\zeta){\rm div}
u_{n}-a\rho_{n}^{\gamma}-\frac{\kappa}{2}\rho_{n}^{2}]\rho_{n}^{\e}=\kappa\rho_{n}^{\e}(-\D)^{-1}{\rm
div}(\rho_{n}(\n\phi*\rho_{n}))-\rho_{n}^{\e}(-\D)^{-1}\p^{2}_{ij}(\rho_{n}
u^{i}{n}u^{j}_{n})\\[2mm]
&\hspace{4,6cm}-\frac{\p}{\p t}\big((\rho_{n})^{\e}(-\D)^{-1}{\rm
div}(\rho_{n}u_{n})\big)+[\frac{\p}{\p t}(\rho_{n})^{\e}]
(-\D)^{-1}{\rm div}(\rho_{n}u_{n}),\\
\end{aligned}
$$
Next we have:
$$
\begin{aligned}
&[(\mu+\zeta){\rm div}
u_{n}-a\rho_{n}^{\gamma}-\frac{\kappa}{2}\rho_{n}^{2}]\rho_{n}^{\e}=\kappa\rho_{n}^{\e}(-\D)^{-1}{\rm
div}(\rho_{n}(\n\phi*\rho_{n}))-\rho_{n}^{\e}(-\D)^{-1}\p^{2}_{ij}(\rho_{n}
u^{i}_{n}u^{j}_{n})\\[2mm]
&\hspace{4,1cm}-\frac{\p}{\p t}[(\rho_{n})^{\e}(-\D)^{-1}{\rm
div}(\rho_{n}u_{n})]- {\rm  div}[ u_{n}(\rho_{n})^{\e}
(-\D)^{-1}{\rm div}(\rho_{n}u_{n})]\\[2mm]
&\hspace{2,5cm}+(\rho_{n})^{\e}u_{n}\cdot\n(-\D)^{-1}{\rm
div}(\rho_{n}u_{n}) +(1-\e)({\rm
div}u_{n})(\rho_{n})^{\e}(-\D)^{-1}{\rm
div}(\rho_{n}u_{n}),\\[2mm]
\end{aligned}
\label{45}
\end{equation}
or finally:
\begin{equation}
\begin{aligned}
&[(\mu+\xi){\rm div}
u_{n}-a\rho_{n}^{\gamma}-\frac{\kappa}{2}\rho_{n}^{2}]\rho_{n}^{\e}=\kappa\rho_{n}^{\e}(-\D)^{-1}{\rm
div}(\rho_{n}(\n\phi*\rho_{n}))\\[2mm]
&\hspace{0,5cm}-\frac{\p}{\p t}[\rho_{n}^{\e}(-\D)^{-1}{\rm
div}(\rho_{n}u_{n})]-{\rm div}[ u_{n}(\rho_{n})^{\e}(-\D)^{-1}{\rm
div}(\rho_{n}u_{n})]\\[2mm]
&\hspace{3cm}+(\rho_{n})^{\e}[u_{n}.\n(-\D)^{-1}{\rm
div}(\rho_{n}u_{n})-(-\D)^{-1}\p^{2}_{ij}(\rho
u^{i}_{n}u^{j}_{n})]\\[2mm]
&\hspace{6,3cm}+(1-\e)({\rm div}u_{n})(\rho_{n})^{\e}(-\D)^{-1}{\rm div}(\rho_{n}u_{n})\;.\\
\end{aligned}
\label{46}
\end{equation}
\\
Now like in Lions \cite{L2} we want to pass to the limit in the
distribution sense in (\ref{46}) in order to estimate
$\overline{{\rm div}u(\rho)^{\e}}$.
\subsubsection*{Passage to the weak limit in (\ref{46})}
In this goal we use the following lemma of P-L Lions in \cite{4L2}
to express the weak limit of the non-linear terms.
\begin{lemme}
\label{4lemme1} Let $\Omega$ a open of $\R^{N}$. Let $(g_{n},h_{n})$
converge weakly to $(g,h)$ in
$L^{p_{1}}(0,T,L^{p_{2}}(\Omega))\times
L^{q_{1}}(0,T,L^{q_{2}}(\Omega))$
where $1\leq p_{1},p_{2},q_{1},q_{2}\leq+\infty$ satisfy,\\
$$\frac{1}{p_{1}}+\frac{1}{q_{1}}=\frac{1}{p_{2}}+\frac{1}{q_{2}}=1\;.$$
We assume in addition that:
\begin{equation}
\frac{\p g^{n}}{\p t}\;\;\mbox{is bounded in}\;\;
L^{1}(0,T,W^{-m,1}(\Omega))\;\;\mbox{for
some}\;\;m\geq0\;\;\mbox{independent of}\;\;n.\hspace{8cm}\label{4101}\\
\end{equation}
and that:
\begin{equation}
\|h^{n}-h^{n}(\cdot,\cdot+\xi)\|_{L^{q_{1}}(0,T,L^{q_{2}}(\Omega))}\rightarrow0\;\;\;\;\mbox{as}\;\;
|\xi|\rightarrow0,\,\mbox{uniformly in n.}\\
\label{4101a}
\end{equation}
Then, $g^{n}h^{n}$ converges to $gh$ (in the sense of distribution
on $\Omega\times(0,T))$.
\end{lemme}
So we use the above lemma to pass to the weak limit in the four following non-linear terms of (\ref{46}):\\
$$
\begin{array}{ll}
T^{1}_{n}=u_{n}\rho_{n}^{\e}(-\D)^{-1}{\rm
div}(\rho_{n}u_{n}),\;\;\;&T^{2}_{n}=\rho_{n}^{\e}(-\D)^{-1}{\rm
div}(\rho_{n}u_{n}),\\[2mm]
T^{3}_{n}=\rho_{n}^{\e}(-\D)^{-1}{\rm
div}(\rho_{n}(\n\phi*\rho_{n})),\;\;\;&
T^{4}_{n}=({\rm div}u_{n})(\rho_{n})^{\e}(-\D)^{-1}{\rm div}(\rho_{n}u_{n})\;.\\
\end{array}
$$
So we choose the different $g_{n}^{i}$ and $h_{n}^{i}$ as follows:
$$
\begin{aligned}
&\mbox{for}\;\;T^{1}_{n}\hspace{1cm}g_{n}^{1}=u_{n}(\rho_{n})^{\e}\hspace{1,3cm}g^{1}=u\overline{\rho^{\e}}\hspace{1,3cm}
h_{n}^{1}=(-\D)^{-1}{\rm
div}(\rho_{n}u_{n})\hspace{3cm}\\
&\mbox{for}\;\;T^{2}_{n}\hspace{1cm}g_{n}^{2}=\rho_{n}^{\e}\hspace{2,15cm}g^{2}=\overline{\rho^{\e}}\hspace{1,5cm}
h_{n}^{2}=(-\D)^{-1}{\rm div}(\rho_{n}u_{n})\\
&\mbox{for}\;\;T^{3}_{n}\hspace{1cm}g_{n}^{3}=\rho_{n}^{\e}\hspace{2,15cm}g^{3}=\overline{\rho^{\e}}\hspace{1,5cm}
h_{n}^{3}=(-\D)^{-1}{\rm
div}(\rho_{n}(\n\phi*\rho_{n}))\\
&\mbox{for}\;\;T^{4}_{n}\hspace{1cm}g_{n}^{4}=({\rm
div}u_{n})(\rho_{n})^{\e}\;\;\hspace{0,3cm}g^{4}=
\overline{{\rm div}u\,\rho^{\e}}\;\;\;\hspace{0,3cm}h_{n}^{4}=(-\D)^{-1}{\rm div}(\rho_{n}u_{n}).\\
\end{aligned}
$$
To show that $u_{n}(\rho_{n})^{\e}$ converges in distribution sense
to $u\overline{\rho^{\e}}$ we apply easily lemma \ref{4lemme1} with
$h_{n}=u_{n}$ and $g_{n}=\rho_{n}^{\e}$.
We now want to examine each term and apply the above lemma to pass
to
the limit in the weak sense.\\
We start with the first term $T^{1}_{n}$. We have that
$\rho_{n}^{\e}u_{n}\in L^{\infty}(L^{q})\cap L^{2}(L^{r})$ with
$\frac{1}{q}=\frac{\e}{2s}+\frac{1}{2}$ and
$\frac{1}{r}=\frac{(N-2)}{2N}+\frac{\e}{s}=\frac{1}{2}-\frac{1}{N}+\frac{\e}{s}$.
In addition the hypothesis (\ref{4101}) is immediately verified (use the momentum equation).\\
We now want to verify the hypothesis $(\ref{4101a})$, so we have
$h^{1}_{n}$ belongs to $L^{\infty}(W_{loc}^{1,q^{'}}(\R^{N}))\cap
L^{2}(W_{loc}^{1,r^{'}}(\R^{N}))$
with $\frac{1}{q^{'}}=\frac{1}{2}+\frac{1}{2s}$ and
$\frac{1}{r^{'}}=\frac{(N-2)}{2N}+\frac{1}{s}=\frac{1}{2}+\frac{1}{s}-\frac{1}{N}$.
This result enables us to verify the hypothesis (\ref{4101a}) by Sobolev embedding.\\
So we can choose (with the notation of the above lemma) $q_{1}=2$
and $q_{2}\in (r^{'},\frac{Nr^{'}}{N-r^{'}})$, $p_{1}=2$, and
$p_{2}=1-\frac{1}{q_{2}}$ which is possible by interpolation. Indeed
we have:
$\frac{1}{r^{'}}+\frac{1}{r}=1-\frac{2}{N}+\frac{1+\e}{s}\leq1$.\\
We proceed in the same way for $T^{2}_{n}$ and $T^{4}_{n}$.
\\
We can similarly examine $T^{3}_{n}$, because $\rho_{n}^{\e}\in
L^{\infty}(L^{\frac{1}{\e}}\cap L^{\frac{s}{\e}})$ and
$\rho_{n}(\n\phi*\rho_{n})\in L^{\infty}(L^{1}\cap
L^{\frac{s}{2}})$, we can choose $p_{2}=\frac{1}{\e}$, we have then
$(\D)^{-1}{\rm div}\rho_{n}(\n\phi*\rho_{n})\in
L^{\infty}(0,T;W^{1,\frac{s}{2}})$ so that we can choose $q_{1}=2$,
$q_{2}\in(1,\frac{N\frac{s}{2}}{N-\frac{s}{2}})$. We can conclude by
interpolation.
\\
\\
Finally we have to study the last non linear following term that we
treat similarly as P-L.Lions in \cite{4L2}:
$$A_{n}=(\rho_{n})^{\e}[u_{n}.\n(-\D)^{-1}{\rm div}(\rho_{n}
 u_{n})-(-\D)^{-1}\p_{ij}^{2}(\rho_{n}
 (u_{i})_{n}(u_{j})_{n})].$$
We can express this term $A_{n}$ as follows:
$$A_{n}=(\rho_{n})^{\e}[u^{j}_{n},R_{ij}](\rho_{n} u^{i}_{n}).$$
where $R_{ij}=(-\D)^{-1}\p_{ij}^{2}$  with $R_{i}$ the classical Riesz transform.\\
Next, we use a result by Coifman, Meyer on this type of commutator
(see \cite{4M}) to take advantage of
 the regularity of $[u^{j}_{n},R_{ij}](\rho_{n} u^{i}_{n})$.
\begin{theorem}
The following map is continuous for any $N\geq 2$:
\begin{equation}
\begin{aligned}
&W^{1,r_{1}}(\R^{N})^{N}\times L^{r_{2}}(\R^{N})\rightarrow W^{1,r_{3}}(\R^{N})^{N}\\
&\hspace{3cm}(a,b)\rightarrow[a_{j},R_{i}R_{j}]b_{i}\\
\end{aligned}
\label{4applicationlineaire}
\end{equation}
with: $\frac{1}{r_{3}}=\frac{1}{r_{1}}+\frac{1}{r_{2}}$.
\end{theorem}
To pass to the weak limit in $A_{n}$ we will use the previous lemma.
We start with the case with $s>3$. This quantity belongs to the
space $L^{1}(W^{1,q})$ provided that $D u_{n}\in L^{2}(L^{2})$ and
$\rho u^{j}\in L^{2}(L^{r})$ where
$\frac{1}{r}=\frac{N-2}{2N}+\frac{1}{s}=\frac{1}{2}-\frac{1}{N}+\frac{1}{s}$
in
which case $\frac{1}{q}=\frac{1}{r}+\frac{1}{2}=1-\frac{1}{N}+\frac{1}{s}\leq1$.\\
for After we can use the above lemma applied to
$h_{n}=[R_{ij},u^{j}_{n}](\rho_{n}u^{i}_{n})$ and
$g_{n}=\rho_{n}^{\e}$. We can show easily in using again lemma
\ref{4lemme1} that $h_{n}$ converges in distribution sense to
$[R_{ij},u_{j}](\rho_{n}u_{i})$.\\
So we can take: $q_{1}=1$, $p_{1}=+\infty$ and $q_{2}\in
(q,\frac{qN}{N-q})$, $p_{2}=1-\frac{1}{q_{2}}$, this one because we
can use interpolation and we can localize as we want limit in distribution sense.\\
In the case where $s\leq 3$, a simple interpolation  argument can be
used to accommodate the general case. It suffices to fix
$L^{r_{2}}(\R^{N})$ in the application (\ref{4applicationlineaire})
and use a result of Riesz-Thorin.
\\
Finally according to the equation (\ref{44}), and after passing to
the limit we get:
\begin{equation}
\begin{aligned}
&[(\mu+\xi)\overline{{\rm div}
u\,\rho^{\e}}-\overline{(a\rho^{\gamma+\e})}-\frac{\kappa}{2}\overline{\rho^{2+\e}}]=\overline{\rho^{\e}}(-\D)^{-1}{\rm
div}(\overline{\rho(\n\phi*\rho)})-\frac{\p}{\p
t}[\overline{\rho^{\e}}(-\D)^{-1}{\rm
div}(\rho u)]\\
&\hspace{2,4cm}-{\rm div}[\overline{\rho^{\e}}u(-\D)^{-1}{\rm
div}(\rho u)]+\overline{\rho^{\e}}[u.\n(-\D)^{-1}{\rm div}(\rho
u)-(-\D)^{-1}\p_{ij}(\rho u_{i}u_{j})]\\
&\hspace{9cm}+(1-\e)\overline{{\rm div} u\,\rho^{\e}}(-\D)^{-1}{\rm
div}(\rho u).
\label{4111a}\\
\end{aligned}
\end{equation}
\subsubsection*{Inequality between the terms $\overline{\rho^{\e}}\,{\rm div}u$ and $\overline{{\rm div}
u\,\rho^{\e}}$} Now we are interested in estimating the term
$\overline{\rho^{\e}}{\rm div}u$ in order to describe the quantity
$\overline{\rho^{\e}}\,{\rm div}u-\overline{{\rm div} u\,\rho^{\e}}$
before considering the quantity $\rho^{\e}\,{\rm
div}u-\overline{{\rm div} u\,\rho^{\e}}$. We pass to the weak limit
directly in (\ref{43}) and we get in using again the lemma
\ref{lemme1}:
\begin{equation}
\begin{aligned}
\frac{\p}{\p t}(-\D)^{-1}{\rm div}(\rho
u)+(-\D)^{-1}\p_{ij}^{2}(\rho u_{i}u_{j})+&[(\mu+\xi){\rm div}
u-\overline{a\rho^{\gamma}}]=\\
&-(-\D)^{-1}{\rm div}(\overline{\rho(\n\phi*\rho)}
+\frac{\kappa}{2}\overline{\rho^{2}}. \label{410a1}
\end{aligned}
\end{equation}
Now we just multiply (\ref{410a1}) with $\overline{\rho^{\e}}$ and
we can see that each term has a distribution
sense.\\
So we get in proceeding in the same way as before:
\begin{equation}
\begin{aligned}
&[(\mu+\xi){\rm div}
u\,\overline{\rho^{\e}}-\overline{(a\rho^{\gamma})}\;\overline{\rho^{\e})}-\frac{\kappa}{2}\overline{\rho^{2}}\overline{\rho^{\e}}]=
\overline{\rho^{\e}}(-\D)^{-1}{\rm
div}(\overline{\rho(\n\phi*\rho))}\\[2mm]
&\hspace{0,3cm}-\overline{\rho^{\e}}\frac{\p}{\p
t}[\rho(-\D)^{-1}{\rm div} (\rho
u)]+\overline{\rho^{\e}}[u.\n(-\D)^{-1}{\rm div}(\rho u)-(-\D)^{-1}\p_{ij}(\rho u_{i}u_{j})]\\[2mm]
&\hspace{3,2cm}-{\rm div}[\rho^{\e}u(-\D)^{-1}{\rm div}(\rho
u)]+(1-\e)\overline{{\rm div}
u(\rho)^{\e}}\,(-\D)^{-1}{\rm div}(\rho u).\\
\end{aligned}
\label{411a}
\end{equation}
In subtracting (\ref{411a}) from (\ref{4111a}), we get:
$$
\begin{aligned}
&(\mu+\xi)\,\overline{{\rm div}
u(\rho)^{\e}}-a\overline{\rho^{\gamma+\e}}-\frac{\kappa}{2}\overline{\rho^{2+\e}}=(\mu+\xi){\rm
div} u\,
\overline{\rho^{\e}}-a\overline{\rho^{\gamma}}\,\overline{\rho^{\e}}-\frac{\kappa}{2}\overline{\rho^{2}}
\overline{\rho^{\e}}\;\;\;\mbox{a.e}\;.\\
\end{aligned}
$$
Next we observe that by convexity:
$$(\overline{\rho^{\gamma+\e}})^{\frac{\e}{\gamma+\e}}\geq(\overline{\rho^{\e}}),\;\;\;
(\overline{\rho^{\gamma+\e}})^{\frac{\gamma}{\gamma+\e}}\geq(\overline{\rho^{\gamma}})\;\;\;\;\mbox{a.e}\;.$$
So we get:
\begin{equation}
\overline{{\rm div}u\,(\rho)^{\e}}\geq{\rm
div}u\,\overline{\rho^{\e}}. \label{420a}
\end{equation}
\subsubsection*{Comparison between $\rho$ and
$\overline{\rho^{\e}}^{\frac{1}{\e}}$} As on since
$(\rho_{n},u_{n})$ are regular solutions we get the equality
(\ref{42}) applied to $B(x)=x^{\e}$. So we get:
\begin{equation}
\frac{\p}{\p t}(\rho_{n})^{\e}+{\rm
div}(u_{n}(\rho_{n})^{\e})=(1-\e){\rm div} u_{n}(\rho_{n})^{\e}.
\label{4e1}
\end{equation}
And after we pass to the weak limit in (\ref{4e1}) and we get:
\begin{equation}
\frac{\p}{\p t}\overline{\rho^{\e}}+{\rm
div}(u\,\overline{\rho^{\e}})=(1-\e)\overline{{\rm div} u\rho^{\e}}.
\label{4102}
\end{equation}
When we pass to the weak limit in combining with (\ref{420a}) we get:\\
\begin{equation}
\frac{\p}{\p t}\overline{(\rho)^{\e}}+{\rm
div}(u\overline{(\rho)^{\e}})\geq(1-\e){\rm div}u
\overline{(\rho)^{\e}}. \label{412}
\end{equation}
Now we wish to conclude about the pointwise convergence of
$\rho_{n}$ in proving that
$(\overline{\rho^{\e}})^{\frac{1}{\e}}=\rho$ and to finish we will
use the following theorem (see \cite{4F} p 34) applied to
$B(x)=x^{\frac{1}{\e}}$ which is convex.
\begin{theorem}
\label{4convexite} Let $(v_{n})_{n\in\mathbb{N}}$ be a sequence of
functions bounded in $L^{1}(\R^{N})$ such that:
$$v_{n}\rightharpoonup v\;\;\;\;\mbox{weakly in}\;\;L^{1}(\R^{N}).$$
Let $\va:\R\longrightarrow(-\infty,+\infty]$ be a upper
semi-continuous strictly concave function such that $\va(v_{n})\in
L^{1}(\R^{N})$ for any $n$, and:
$$\va(v_{n})\rightharpoonup\overline{\va(v)}\;\;\;\mbox{weakly in}\;\;L^{1}(\R^{N}).$$
Then:
$$v_{n}(y)\rightarrow v(y)\;\;\mbox{a.e.}$$
extracting a subsequence as the case may be.
\end{theorem}
Now we want to use a type of Diperna-Lions theorem on inequality
(\ref{412}). Our goal is to renormalize this inequality with the
function $B(x)=x^{\frac{1}{\e}}$ so that one can compare $\rho$ and
$\overline{\rho^{\e}}^{\frac{1}{\e}}$. Although (\ref{412}) doesn't
correspond exactly to the mass equation, we can use the same
technics to renormalize the solution provided that $\rho\in
L^{\infty}(L^{2})$ which is the case. In our case it is very
important that $\rho\in L^{\infty}(L^{2})$, indeed it avoids to have
supplementary conditions on the index $\gamma$ like for the
compressible Navier-Stokes system in \cite{4L2}. We recall
 of Diperna-Lions theorem on renormalized solution
for the mass equation.
\begin{theorem}
Suppose that $\rho\in L^{\infty}(L ^{2})$ and $\beta\in C^{1}([0,\infty);\R)$.\\
We have then:
$$\frac{\p \beta(\rho)}{\p t}+{\rm div}(\beta(\rho)\,u)=(\beta(\rho)-\rho\beta^{'}(\rho)){\rm div}u$$
in distribution sense.
\end{theorem}
We now want to adapt this theorem for our equation (\ref{12}) with
$\beta(x)=x^{\frac{1}{\e}}$ , so we may regularize by
$\omega_{\alpha}$ (with
$\omega_{\alpha}=\frac{1}{\alpha^{N}}\omega(\frac{\cdot}{\alpha})$
where $\omega\in C_{0}^{\infty}(\R^{N})$, $\mbox{supp}\;\omega\in
B_{1}$ and $\int\omega dx=1$) and find for all $\beta\in
C^{\infty}_{0}([0,+\infty))$:
$$
\begin{aligned}
&\frac{\p}{\p t}(\overline{\rho^{\e}}*\omega_{\alpha})+{\rm
div}[u\,\overline{\rho^{\e}}*\omega_{\alpha}]\geq (1-\e){\rm div}u\,
\overline{\rho^{\e}}*\omega_{\alpha}+R_{\alpha}\\
\end{aligned}
$$
where we have:
$$R_{\alpha}={\rm
div}[u\,\overline{\rho^{\e}}*\omega_{\alpha}]-{\rm
div}(u\,\overline{\rho^{\e}})*\omega_{\alpha}+(1-\e)[{\rm div}u\,
\overline{\rho^{\e}}]*\omega_{\alpha}-(1-\e){\rm div}u\,
\overline{\rho^{\e}}*\omega_{\alpha}$$
$$
\begin{aligned}
&\frac{\p}{\p t}(\beta(\overline{\rho^{\e}}*\omega_{\alpha}))+{\rm
div}[u\,\beta(\overline{\rho^{\e}}*\omega_{\alpha})]\geq (1-\e){\rm
div}u\,
\overline{\rho^{\e}}*\omega_{\alpha}\,\beta^{'}(\overline{\rho^{\e}}*\omega_{\alpha})\\[2mm]
&\hspace{3.6cm}+({\rm
div}u)[\beta(\overline{\rho^{\e}}*\omega_{\alpha})-\overline{\rho^{\e}}*\omega_{\alpha}\beta^{'}
(\overline{\rho^{\e}}*\omega_{\alpha})]+R_{\alpha}\beta^{'}(\overline{\rho^{\e}}*\omega_{\alpha})\\[2.5mm]
&\hspace{5cm}=-\e({\rm
div}u)\overline{(\rho)^{\e}}\beta^{'}(\overline{\rho^{\e}})+({\rm
div}u)\,\beta(\overline{\rho^{\e}})+R_{\alpha}\beta^{'}(\overline{\rho^{\e}}*\omega_{\alpha})\;.\\
\end{aligned}
$$
After we pass to the limit when $\alpha\rightarrow0$ and we see that
$R_{\alpha}$ tends to $0$ in using lemma on regularization in
\cite{4L1} p 43. This looks like a rather innocent manipulation but
it's at this point that we require to control $\rho$ in
$L^{2}(0,T;\R^{N})$. And in our case we don't need to impose
$\gamma>\frac{N}{2}$ for $N=2,3$. Hence:
$$
\begin{aligned}
\frac{\p}{\p t}(\beta(\overline{(\rho)^{\e}}))+{\rm
div}[u\,\beta(\overline{(\rho)^{\e}})]\geq
&-\e({\rm
div}u)\overline{\rho^{\e}}\beta^{'}(\overline{\rho^{\e}})+({\rm
div}u)\beta(\overline{\rho^{\e}}).\\
\end{aligned}
$$
\\
We then choose $\beta=(\Psi_{M})^{\frac{1}{\e}}$ where
$\Psi_{M}=M\Psi(\frac{\cdot}{M})$, $M\geq1,\,\Psi\in
C^{\infty}_{0}([0,+\infty)),\,\Psi(x)=x$ on $[0,1]$,
$\mbox{supp}\Psi\subset[0,2]$, and we obtain:\\
$$
\begin{aligned}
&\frac{\p}{\p t}(\Psi_{M}(\overline{\rho^{\e}})^{\frac{1}{\e}})+{\rm
div}[u\,\Psi_{M}(\overline{\rho^{\e}})^{\frac{1}{\e}}]\\[2mm]
&\;\;\;\;\;\;\;\;\;\;\;\;\;\;\;\;\;\;\;\;\;\;\;\;\;\;\;\;\geq({\rm
div}u)\Psi_{M}(\overline{\rho^{\e}})^{\frac{1}{\e}-1}\Psi_{M}^{'}(\overline{\rho^{\e}})\overline{\rho^{\e}}+({\rm
div}u)\Psi_{M}(\overline{\rho^{\e}})^{\frac{1}{\e}}\\[2mm]
&\;\;\;\;\;\;\;\;\;\;\;\;\;\;\;\;\;\;\;\;\;\;\;\;\;\;\;\;\geq{\rm
div}u\Psi_{M}(\overline{\rho^{\e}})^{\frac{1}{\e}-1}
[\Psi_{M}(\overline{\rho^{\e}})^{\frac{1}{\e}}-\Psi_{M}^{'}(\overline{\rho^{\e}})\overline{\rho^{\e}}]
1_{(\overline{\rho^{\e}}>M)}\\[2mm]
&\;\;\;\;\;\;\;\;\;\;\;\;\;\;\;\;\;\;\;\;\;\;\;\;\;\;\;\;\geq
-C_{0}|{\rm div}u|M^{\frac{1}{\e}}1_{(\overline{\rho^{\e}}>M)}.\\
\end{aligned}
$$
\\
where $C_{0}=\sup\{|\Psi(x)|^{\frac{1}{\e}-1}|\Psi(x)-x\Psi^{'}(x)|,
\;x\in[0,+\infty)\}$.\\
Now we claim that:
\begin{equation}
\frac{\p}{\p t}(\overline{(\rho)^{\e}})^{\frac{1}{\e}})+{\rm
div}(u\,\overline{(\rho)^{\e}})^{\frac{1}{\e}})\geq0.\\
\label{413}
\end{equation}
For proving that, we notice that by convexity
$\overline{(\rho)^{\e}}^{\frac{1}{\e}}\leq\rho$, so we get :
$$\||{\rm
div}u|M^{\frac{1}{\e}}1_{\overline{(\rho)^{\e}}>M}\|_{L^{1}_{T}(L^{1}(\R^{N})}\leq\|{\rm
div}u\|_{L^{2}_{T}(L^{2}(\R^{N})}\|\rho\,1_{\rho>M^{\frac{1}{\e}}}\|_{L^{2}_{T}(L^{2}(\R^{N})}\rightarrow0\;\;
\mbox{as}\;\;M\rightarrow +\infty.$$ We have concluded by dominated
convergence.
We have concluded by dominated convergence.
\\
At this stage we subtract the mass equation to (\ref{413}) and we
get in setting $r=\rho-\overline{(\rho)^{\e}}^{\frac{1}{\e}}$:\\
\begin{equation}
\frac{\p}{\p t}(r)+{\rm
div}(ur)\leq0.\\
\label{414}
\end{equation}
We now want to integrate and to use the fact that $r\geq0$ to get
that $r=0$ a.a. To justify the integration we  test our inequality
against a cut-off function of the form $\va(\frac{\cdot}{R})$ where
$\va\in C^{\infty}_{0}(\R^{N})$, $\va=1$ on $B(0,1)$,
$\mbox{Supp}\va\subset B(0,2)$ and $R>1$. We test the equation
$(\ref{414})$ against $\va_{R}$ and we get:
\begin{equation}
\int_{[0,T]\times\R^{N}}\frac{\p}{\p t}[r(t,x)]\va_{R}(x)-u(t,x)r(t,x)\frac{1}{R}\n\va(\frac{x}{R})dt\,dx\leq0.\\
\label{415}
\end{equation}
Next we notice that:
$$
\begin{aligned}
&|\int_{[0,T]\times\R^{N}}u(t,x)r(t,x)\frac{1}{R}\n\va(\frac{x}{R})dt\,dx|\leq\|u\|_{L^{1}(0,T;L^{\frac{2N}{N-2}}(\R^{N}))}
\|r\|_{L^{1}(0,T;L^{\frac{2N}{N+2}}(\R^{N}))}\\
&\hspace{11cm}*\frac{1}{R}\|\n\va\|_{L^{\infty}(\R^{N})}.\\
\end{aligned}
$$
It implies that:
$$\int_{[0,T]\times\R^{N}}u(t,x)r(t,x)\frac{1}{R}\n\va(\frac{x}{R})dt\,dx\rightarrow0\;\;\;\mbox{as}
\;\;R\rightarrow+\infty.$$ We have then:
$$\int_{[0,T]\times\R^{N}}\frac{\p}{\p t}r(t,x)\va_{R}(x)dt\,dx=\int_{\R^{N}}r(T,x)\va_{R}(T,x)dx-
\int_{\R^{N}}r(0,x)\va_{R}(0,x)dx.$$ Now we want to verify that
$r(0,\cdot)=0$ to conclude because. Indeed we will obtain that:
$$\lim_{R\rightarrow+\infty}\int_{\R^{N}}r(T,x)\va_{R}(T,x)dx\rightarrow\int_{\R^{N}}r(T,x)dx\leq0\;\;\;\mbox{and}
\;\;\;r\geq0.$$
then $r=0$.\\
We know that $\rho_{n}$ is uniformly bounded in
$L^{\infty}(L^{1}\cap L^{s}(\R^{N}))$, then $\rho_{n}^{\e}$ is
relatively compact in $C([0,T];L^{p}-w)$ with $1<p<s$ (where
$L^{p}-w$ denote the space $L^{p}$ endowed with weak topology).
Moreover $(\rho_{0}^{\e})_{n}$ converges to $\rho_{0}^{\e}$, we
deduce then $r(0)=0$ a.a.
\\
Now as $r=0$ we conclude in using the theorem \ref{4convexite}
$\rho_{n}$ converges a.a to $\rho$ and that $\rho_{n}$ converges to
$\rho$ in $L^{p}([0,T]\times B_{R})$ for all $p\in[1,q)$ and in
$L^{p_{1}}(0,T,L^{p_{2}}(B_{R}))$ for all $p_{1}\in[1,+\infty)$,
$p_{2}\in[1,s)$ and for all $R\in(0,+\infty)$.
\subsubsection*{Conclusion}
We wish now conclude and get the convergence of our theorem in the
total space.
\\
We aim at proving here the convergence of $\rho_{n}$ in
$C([0,T],L^{p}(\R^{N}))\cap L^{q^{'}}(\R^{N}\times (0,T))$ for all
$1\leq p<s,\,1\leq q^{'}<q$. We have just to show the convergence of
$\rho_{n}$ to $\rho$ in $C([0,T],L^{1}(\R^{N}))$. To this end, we
introduce $d_{n}=\sqrt{\rho_{n}}$ which clearly converges to
$\sqrt{\rho}$ in $L^{2p_{1}}(0,T,L^{2p_{2}}(B_{R}))\cap
L^{2p}(B_{R}\times (0,T))$
to $d=\sqrt{\rho}$ for all $R\in(0,+\infty)$.\\
We next remark that $\rho\in C([0,T],L^{1}(\R^{N}))$ and thus $d\in
C([0,T],L^{2}(\R^{N})$. Indeed, using once more the regularization
lemma in \cite{4L1} we obtain the existence of a bounded
$\rho_{\al}\in C([0,T],L^{1}(\R^{N}))$ smooth in $x$ for all $t$
satisfying:
$$\frac{\p \rho_{\al}}{\p t}+{\rm div}(u\rho_{\al})=r_{\al}\;\;\;\mbox{in}\;\;L^{1}((0,T)\times\R^{N})\;\;\mbox{as}
\;\;\al \rightarrow 0_{+}.$$
with $r_{\al}={\rm div}(u\rho_{\al})-{\rm div}(\rho u)*w_{\al}$ (where $w$ is defined as in the previous part).\\
$\rho_{\al}\rightarrow\rho$ in $L^{1}(\R^{N}\times(0,T))$,
$\rho_{\al}/ _{t=0}\rightarrow\rho/ _{t=0}$ in $L^{1}(\R^{N})$
as $\al\rightarrow0_{+}.$\\
From these facts, it is straightforward to deduce that:
$$\frac{\p}{\p t}|\rho_{\al}-\rho_{\eta}|+{\rm div}(u|\rho_{\al}-\rho_{\eta}|)=r_{\al}-r_{\eta}$$
and thus:
$$\sup_{[0,T]}\int_{\R^{N}}|\rho_{\al}-\rho_{\eta}|dx=\int_{0}^{T}\int_{\R^{N}}|r_{\al}-r_{\eta}|dx.$$
Since $\rho\in C([0,T],L^{p}(B_{R})-w)$ (for all $R\in(0,+\infty)$,
$1<p<s$),
we may then deduce that $\rho_{\al}$ converge to $\rho$ in $C([0,T],L^{1}(\R^{N}))$.\\
Next, we observe that we can justify as we did above that $d_{n}$ and $d$ satisfy:\\
$$\frac{\p d_{n}}{\p t}+{\rm div}(u_{n}d_{n})=\frac{1}{2}d_{n}{\rm div}(u_{n}),$$
$$\frac{\p d}{\p t}+{\rm div}(ud)=\frac{1}{2}d{\rm div}(u).$$
Therefore once more, $d_{n}$ converges to $d$ in $C([0,T],L^{2}(\R^{N})-w)$.\\
Thus in order to conclude, we just have to show that whenever
$t_{n}\in[0,T]$, $t_{n}\rightarrow t$, then $d_{n}(t_{n})\rightarrow
d(t)$ in $L^{2}(\R^{N})$ or equivalently that
$\int_{\R^{N}}d_{n}(t_{n})^{2}dx=\int_{\R^{N}}\rho_{n}(t_{n})dx\rightarrow_{n}\int_{\R^{N}}d(t)^{2}dx=
\int_{\R^{N}}\rho(t)dx$. This is the case since we deduce from the
mass equation, integrating this equation over $\R^{N}$ and
justifying the integration exactly like previously that:
$$\int_{\R^{n}}\rho_{n}(t_{n})dx=\int_{\R^{n}}(\rho_{0})_{n}dx\rightarrow_{n}\int_{\R^{n}}\rho_{0}dx=
\int_{\R^{n}}\rho(t)dx.$$ We then conclude by uniform continuity
that $\|\rho_{n}(t_{n})-\rho_{n}(t)\|_{L^{1}}$ tends to $0$.
\subsubsection*{Case $N=2$}
First of all, the main difficulty is the fact that we no longer have
global $L^{p}$ bounds on $u_{n}$. That's why most of the proof is in
fact local and we know that $u_{n}$ is bounded in
$L^{2}(0,T;L^{p}(B_{R}))$ for all $p\in[1,+\infty)$,
$R\in(0,+\infty)$.\\
As we need to localize the argument, we get the following limit:
$$
\begin{aligned}
&((\mu+\xi){\rm
div}u_{n}-a\rho_{n}^{\gamma}-\frac{1}{2}\rho_{n}^{2})\rho_{n}^{\e}\rightharpoonup_{n}
((\mu+\xi){\rm
div}u-a\overline{\rho^{\gamma}}-\frac{1}{2}\overline{\rho^{2}})\,\overline{\rho^{\e}}\;\;
\mbox{in}\;\;{\cal D}^{'}(\R^{N}\times[0,T])\\
\end{aligned}
$$
let $\varphi\in C_{0}^{\infty}(\R^{N})$, $0\leq\varphi\leq1$,
$\mbox{supp}\varphi\subset K$ for an arbitrary compact set
$K\in\R^{N}$.
We apply the operator $(\D)^{-1}{\rm div}$ to the momentum equation
that we have localized and we pass directly to the weak limit:
\begin{equation}
\begin{aligned}
\frac{\p}{\p t}(-\D)^{-1}{\rm div}(\va\rho u)+R_{ij}(\va\rho
u_{i}u_{j})+&[(\mu+\xi)\rm div
u\,\va-\overline{a\rho^{\gamma}}\va]\\
&=\kappa(-\D)^{-1}{\rm div}(\va\overline{\rho(\n\phi*\rho)}
&-\frac{\kappa}{2}\va\overline{\rho^{2}}))+(-\D)^{-1}\overline{R}.
\label{410}
\end{aligned}
\end{equation}
with:
$$
\begin{aligned}
&\overline{R}=\p_{i}\va\p_{j}(\rho u_{i}u_{j})+(\p_{ij}\va)\rho
u_{i}u_{j}-(\mu+\xi)\D\va{\rm div}u-(2\mu+\xi)\n\va\cdot
\n{\rm div}u+\mu\D u\cdot\n\va\\
&+\D\va
a\overline{\rho^{\gamma}}+a\n\va\cdot\n\overline{\rho^{\gamma}}
+\frac{\kappa}{2}\D\va
\overline{\rho^{2}}+\frac{\kappa}{2}\n\va\cdot\n\overline{\rho^{2}}
.\\
\end{aligned}
$$
Now we  multiply (\ref{410}) with $\overline{\rho^{\e}}$ and we verify that each term has a sense.\\
So we get in proceeding in the same way as before, we can verify
that $(\rho_{n})^{\e}(-\D)^{-1}R_{n}$ converges in distribution
sense to $\overline{\rho^{\e}}(-\D)^{-1}\overline{R}$ for small
enough $\e$.
We get as in the previous case for $N\geq3$:
$$
\begin{aligned}
&\va[(\mu+\xi)\,\overline{{\rm div}
u(\rho)^{\e}}-a\overline{\rho^{\gamma+\e}}-\frac{\kappa}{2}\overline{\rho^{2+\e}}]=\va[(\mu+\xi)({\rm
div} u)\va\,
\overline{\rho^{\e}}-a\overline{\va\rho^{\gamma}}\,\overline{\rho^{\e}}-\frac{\kappa}{2}\overline{\va\rho^{2}}
\overline{\rho^{\e}}]\;\;\;\mbox{a.e}\;.\\
\end{aligned}
$$
We then deduce the following inequalities as in the previous proof:
$$\frac{d}{dt}(\overline{\rho^{\e}})^{\frac{1}{\e}}+{\rm
div}(u(\overline{\rho^{\e}})^{\frac{1}{\e}})\geq0\;\;\;\mbox{in}\;\;{\cal
D}^{'}((0,T)\times\R^{N}).$$ We see that the only point left to
check is the justification of the integration over $\R^{2}$ of terms
like ${\rm div}(\overline{\rho^{\e}}^{\frac{1}{\e}}u)$ or ${\rm
div}(\rho u)$ and more precisely that the integral vanishes. This is
in fact straightforward provided we use the bounds on $\rho\in
L^{\infty}(L^{1}(\R^{N}))$ and $\rho |u|^{2}\in
L^{\infty}(L^{1}(\R^{N}))$ and so $\rho u\in
L^{\infty}(L^{1}(\R^{N}))$. Then, letting $\va\in
C_{0}^{\infty}(\R^{2})$, $0\leq\va\leq1$, $\va=1$ on $B(0,1)$ and
$\va=0$ on $^{c}B(0,2)$. We set
$\va_{R}(\cdot)=\va(\frac{\cdot}{R})$ for $R\geq1$, we have
similarly as in the previous case:
$$\begin{aligned}
&|\int^{T}_{0}dt\int_{\R^{2}}ru\cdot\n\va_{R}(x)dx|\leq\|\n\va\|_{L^{\infty}(\R^{2})}\frac{1}{R}\|\rho
u\|_{L^{1}(
0,T)\times\R^{2})}\rightarrow0\;\;\;\mbox{as}\;\;R\rightarrow+\infty\\
\end{aligned}
$$
We can then conclude as in the previous proof.\\
\hfill{$\Box$}
\subsubsection*{Proof of the convergence assertion on $\rho_{n}u_{n}$}
We now want to show the convergence of $\rho_{n}u_{n}$ to have
informations on strong convergence of $u_{n}$ modulo the vacuum. We
recall in this part some classical inequalities to get the
convergence of $\rho_{n}u_{n}$, for more details see Lions in
$\cite{4L2}$ We use once more a mollifier
$k_{\al}=\frac{1}{\al^{N}}k(\frac{\cdot}{\al})$ where $k\in
C_{0}^{+\infty}(\R^{N})$ and we let $g_{\al}=g*k_{\al}$ for an
arbitrary function $g$. We first observe that we have for all
$\frac{N}{2}<p<s$:
$$
\begin{aligned}
|\big((\rho_{n}u_{n})_{\alpha}-\rho_{n}u_{n}\big)
(x)|=\big|\g[\rho_{n}(t,y)-\rho_{n}(t,x)]&u_{n}(t,y)k_{\alpha}(x-y)dy\\
&+\rho_{n}(t,x)\big((u_{n})_{\al}-u_{n}\big)(t,x)\big|\\
\end{aligned}
$$
We have in using H\"older inequalities with the measure
$k_{\alpha}(x-y)dy$:
$$
\begin{aligned}
&|\big((\rho_{n}u_{n})_{\alpha}-\rho_{n}u_{n}\big)
(x)|\leq\big[\g|\rho_{n}(t,y)-\rho_{n}(t,x)|^{p}k_{\alpha}(x-y)dy\big]^{\frac{1}{p}}\,\big(
|u_{n}|^{\frac{p}{p-1}}\big)_{\alpha}^{\frac{p-1}{p}}\\
&\hspace{10cm}+\rho_{n}|(u_{n})_{\al}-u_{n}|(t,x).\\
\end{aligned}
$$
Hence for all $t\geq 0$
$$
\begin{aligned}
&\g|\big((\rho_{n}u_{n})_{\alpha}-\rho_{n}u_{n}\big)
(x)|dx\leq\big[\g
dx\g|\rho_{n}(t,y)-\rho_{n}(t,x)|^{p}k_{\alpha}(x-y)dy
\big]^{\frac{1}{p}}\\
&\hspace{6cm}
*\|\big(
|u_{n}|^{\frac{p}{p-1}}\big)_{\alpha}\|_{L^{1}}^{\frac{p-1}{p}}+\|\rho_{n}\|_{L^{p}}\|(u_{n})_{\al}-u_{n}
\|_{L^{\frac{p-1}{p}}},\\[3mm]
&\hspace{3cm}\leq\big[\sup_{|z|\leq\al}\|\rho_{n}(\cdot+z)-\rho_{n}\|_{L^{p}}\big]\|u_{n}\|_{L^{\frac{p-1}{p}}}
+\|\rho_{n}\|_{L^{p}}\|(u_{n})_{\al}-u_{n} \|_{L^{\frac{p-1}{p}}}.\\
\end{aligned}
$$
Next if we choose $p>\frac{2N}{N+2}$, $\frac{p}{p-1}<\frac{2N}{N-2}$
then $\|(u_{n})_{\al}-u_{n} \|_{L^{2}(0,T;L^{\frac{p-1}{p}})}$
converge to $0$ as $\al$ goes to $0_{+}$ uniformly in $n$. In
addition, the convergence on $\rho_{n}$ assure that
$\sup_{|z|\leq\al}\|\rho_{n}(\cdot+z)-\rho_{n}\|_{L^{p}}$ converge
to $0$ as $\al$ goes to $0_{+}$ uniformly in $n$. Therefore in
conclusion, $(\rho_{n}u_{n})_{\alpha}-\rho_{n}u_{n}$ converge to $0$
in $L^{2}(0,T;L^{1})$ as $\al$ goes to $0_{+}$ uniformly in $n$.\\
Next $(\rho_{n}u_{n})_{\alpha}$ is smooth in $x$, uniformly in $n$
and in $t\in[0,T]$. Therefore, remarking that $\frac{\p}{\p
t}(\rho_{n}u_{n})_{\alpha}$ is bounded in a $L^{2}(0,T;H^{m})$ for
any $m\geq0$, we deduce that $(\rho_{n}u_{n})_{\alpha}$ converge to
$(\rho u)_{\alpha}$ as $n$ goes to $+\infty$ in
$L^{1}((0,T)\times\R^{N})$ for each $\al$. Then using the bound on
$\rho_{n}u_{n}$ in $L^{\infty}(L^{\frac{2s}{s+1}})$, we deduce that
$\rho_{n}u_{n}$ converge to $\rho u$ in $L^{1}((0,T)\times\R^{N})$
and we can
conclude by interpolation.\\
The last convergence is a consequence of the strong convergence of
$\rho_{n}$ and $\rho_{n}u_{n}$.
\section{Existence of weak solution with  general pressure}
In the sequel we concentrate us only on the cases $N=2,3$. We now
want to extend our previous result to more general and physical
pressure laws. In particular we are now interested by two cases, the
first one concerns monotonous pressure law (close in a certain sense
that we will precise to  $\rho^{\gamma}$ pressure) , the second one
is the case of a slightly modified Van der Waals pressure.\\
The technics of proof will be very similar to the previous proof,
only technical points change.
\subsection{Monotonous pressure}
In this section, we shall investigate an extension of the preceding
results  to the case of a general monotonous pressure $P$, i.e $P$
is
assumed to be a $C^{1}$ non-decreasing function on $[0,+\infty)$ vanishing at $0$.\\
We want here to mention in the general situation our new energy
inequality, we recall the inequality (\ref{4inegaliteenergie}):
$$
\begin{aligned}
&\int_{\R^{N}}(\frac{1}{2}\rho
|u|^{2}+\Pi(\rho)+E_{global}[\rho(.,t)])(x)dx(t)+\int_{0}^{t}\int_{\R^{N}}(\mu
D(u):D(u)\\
&\hspace{2cm}+(\lambda+\mu)|{\rm div} u|^{2})dx
\leq\int_{\R^{N}}\big(\frac{|m_{0}|^{2}}{2\rho}+\Pi(\rho_{0})+E_{global}[\rho_{0}]\big)dx.
\end{aligned}
$$
where we define $\Pi$ by:
$\frac{\p}{\p t}(\frac{\Pi(t)}{t})=\frac{P(t)}{t^{2}}$ for all $t>0$.\\
There are two cases worth considering: first of all if $P(t)$ is
such that $\int^{1}_{0}\frac{P(s)}{s^{2}}ds<+\infty$ then we can
choose $\Pi(\rho)=\rho\int^{\rho}_{0}\frac{P(s)}{s^{2}}ds$.\\
In the other case, i.e $\lim_{t\rightarrow0}\frac{P(t)}{t}=c>0$, we
can choose $\Pi(\rho)=\rho\int^{1}_{\rho}\frac{p(s)}{s^{2}}ds$ and
$\Pi$ behaves like $\rho\log\rho$ as $\rho$ goes to $0$.
\\
We now consider a sequence of solutions $(\rho_{n},u_{n})$ and we
make the same assumptions on this sequence as in the previous
section except that we need to modify the assumptions on $\rho_{n}$.
We assume always that $\rho_{n}$ is bounded in
$C([0,T],L^{1}(\R^{N}))$, $P(\rho_{n})$ is bounded in
$L^{\infty}(0,T,L^{1}(\R^{N}))$.
$$
\begin{aligned}
&.\;(\rho_{n})_{n\geq1}\;\;\mbox{is bounded in}\;\;
L^{\infty}(0,T,L^{2}),\hspace{10cm}
\end{aligned}
$$
and we also assume that we have:
$$(\rho_{n}^{\e}P(\rho_{n}))_{n\geq1}\;\;\mbox{is bounded
in}\;\;L^{1}(K\times(0,T))$$ for some $\e>0$, where $K$ is an
arbitrary compact set included in $\R^{N}$.
\\
\begin{theorem}
Let the assumptions of theorem \ref{4principal} be satisfied with in
addition
P a monotone pressure.\\
Then there exists a renormalized finite energy weak solution to
problem $(NSK)$ in the sense of definitions \ref{4defrenormal} and
\ref{4defsolutionsfaibles}. Moreover $P(\rho_{n})$ converges to
$P(\rho)$ in $L^{1}(K\times(0,T))$ for any compact set $K$.
\label{4T32}
\end{theorem}
{\bf Proof:}\\
\\
In this situation, the proof of theorem \ref{4T32} is based on the
same compactness argument as in the theorem \ref{4principal}. In
particular, there is essentially one observation which allows us to
adapt the proof of theorem  \ref{4principal}. Namely we stills
obtain the following identity for the effective viscous flux:
$$\beta(\rho_{n})((\mu+\xi){\rm
div}u_{n}-P(\rho_{n})-\frac{\kappa}{2}\rho_{n}^{2
})\rightharpoonup_{n}\overline{\beta(\rho)}((\mu+\xi){\rm
div}u-\overline{P(\rho)}-\frac{\kappa}{2}\overline{\rho}^{2}),$$
with $\beta(\rho)=\rho^{\e}$ for $\e$ small enough. We have then:
\begin{equation}
(\mu+\xi)\overline{\rho^{\e}{\rm
div}u}-\overline{P(\rho)\rho^{\e}}-\frac{\kappa}{2}\overline{\rho^{2+\e
}})=(\mu+\xi)\overline{{\rm
div}u}\overline{\rho^{\e}}-\overline{P(\rho)}\overline{\rho^{\e}}-\frac{\kappa}{2}\overline{\rho^{2}}\overline{\rho^{\e}}),
\label{4effectifmono}
\end{equation}
Now we can recall a lemma coming from P.-L.Lions in \cite{4L2}:
\begin{lemme}
Let $p_{1}$, $p_{2}\in C([0,\infty))$ be non-decreasing functions.
We assume that $p_{1}(\rho_{n})$, $p_{2}(\rho_{n})$ are relatively
weakly compact in $L^{1}(K\times(0,T))$ for any compact set
$K\subset\R^{N}$. Then, we have:
$$\overline{p_{1}(\rho)p_{2}(\rho)}\geq\overline{p_{1}(\rho)}\overline{p_{2}(\rho)}\;\;\mbox{a.e}.$$
\label{4lemme2}
\end{lemme}
We get finally as in the proof of theorem \ref{4principal} in using
lemma \ref{4lemme2}:
$$\overline{{\rm div}u\,\rho^{\e}}\geq{\rm div}u\overline{\rho^{\e}},$$
All remaining argumentation of the proof of theorem \ref{4principal}
can be performed to conclude. \hfill{$\Box$}
\subsection{Pressure of Van der Waals type}
In this section we are interested by pressure of type Van der Waals
which consequently are not necessarily non-decreasing. That's why in
the
following proof we will proceed slightly differently.\\
So we consider the pressure law:
$$
\begin{aligned}
P(\rho)&=\frac{RT_{*}\rho}{b-\rho}-a\rho^{2}\;\;\mbox{for}\;\;\rho\leq
b-\theta\;\;\;\mbox{for some small}\;\;\theta>0\\[2mm]
\end{aligned}
$$
and we extend the function $P$ to be strictly increasing function
on $\rho\geq b-\theta$.\\
We have then:
\begin{enumerate}
\item $P^{'}$ is bounded from below, that is:
$$P^{'}(\rho)\geq-\bar{\rho}\;\;\;\mbox{for all}\;\;\rho>0.$$
\item $P$ is a strictly increasing function for $\rho$ large enough.
\end{enumerate}
Under the above conditions, it is easy to see that the pressure can
be written as:
$$P(\rho)=P_{1}(\rho)-P_{2}(\rho).$$
with $P_{1}$ a non-decreasing function of $\rho$, and
$$P_{2}\in C^{2}[0,+\infty),\;\;P_{2}\geq0,\;\;P_{2}=0\;\;\mbox{for}\;\;\rho\geq\overline{\rho}.$$
\begin{rem}
The a priori energy estimate give us the bound of $\rho$ in
$L^{\infty}(L^{2})$, then we have:
$$|\{(t,x)\in(0,T)\times\R^{N}/|\rho(t,x)|>b\}|\leq\frac{T\|\rho\|^{2}_{L^{\infty}(L^{2}(\R^{N}))}}{b^{2}}.$$
Then we can controll in measure the set where $P$ is different from
the Van der Waals measure, and it's a set of finite measure.
\end{rem}
We obtain the following theorem.
\begin{theorem}
If in addition to the above assumptions, we assume that
$\rho^{0}_{n}$ converges in $L^{1}(\R^{N})$ to $\rho_{0}$ then
$(\rho,u)$ is a weak solution of the system $(NSK)$
satisfying the initial condition and we have:\\
$$
\begin{aligned}
\rho_{n}\rightarrow\rho\;\;\mbox{in}\;\;C([0,T],L^{p}(\R^{N})\cap
L^{r}((0,T)\times\R^{N})\;\;
&\mbox{for all}\;\;1\leq p<2,\,1\leq r<1+\frac{4}{N},\\
\end{aligned}
$$
\label{4theorem3.8}
\end{theorem}
{\bf Proof: }\\
\\
Most of the proof of theorem \ref{4theorem3.8} is similar as theorem
\ref{4principal}. We will use a approximated sequel $T_{k}$
(introduced by Feireisl in \cite{F}) of $\rho$ by some concave
bounded function.
\begin{definition}
We define the function $T\in C^{\infty}(\R^{N})$ as follows:
$$
\begin{aligned}
&T(z)=z\;\;\;\mbox{for}\;\;z\in[0,1],\\
&T(z)\;\;\mbox{concave on}\;\;[0,+\infty),\\
&T(z)=2\;\;\;\mbox{for}\;\;z\geq3,\\
&T(z)=-T(-z)\;\;\;\mbox{for}\;\;z\in(-\infty,0],\\
\end{aligned}
$$
And $T_{k}$ is the cut-off function:
$$T_{k}(z)=kT(\frac{z}{k}).$$
\end{definition}
In following the proof of theorem \ref{4principal}, we get:
$$
\begin{aligned}
&\frac{\p}{\p t}(\overline{L_{k}(\rho)}-L_{k}(\rho))+{\rm
div}((\overline{T_{k}(\rho)}-T_{k}(\rho))u)+\overline{T_{k}(\rho){\rm
div}u}-\overline{T_{k}(\rho)}{\rm
div}u\\
&\hspace{7cm}=(T_{k}(\rho)-\overline{T_{k}(\rho)}){\rm
div}u\;\;\;\mbox{in}\;\;{\cal D}^{'}((0,T)\times\R^{N}),\\
\end{aligned}
$$
where $L_{k}(\rho)=\rho\log(\rho)$ for $0\leq\rho\leq k$, and $0\leq L_{k}(\rho)\leq \rho\log(\rho)$ otherwise.\\
So we get in integrating in time on $[t_{1},t_{2}]$:
$$
\begin{aligned}
&\int_{\R^{N}}\big(\overline{L_{k}(\rho)}-L_{k}(\rho)\big)(t_{2})dx-\int_{\R^{N}}\big(\overline{L_{k}(\rho)}-L_{k}(\rho)\big)(t_{1})dx\\
&\hspace{1,5cm}+\int^{t_{2}}_{t_{1}}\int_{\R^{N}}\overline{T_{k}(\rho){\rm
div}u}-\overline{T_{k}(\rho)}{\rm div}u\,dxdt=\int^{t_{2}}_{t_{1}}
\int_{\R^{N}}(T_{k}(\rho)-\overline{T_{k}(\rho)}){\rm
div}u\,dxdt.\\
\end{aligned}
$$
We can show that:
$$\|T_{k}(\rho)-\overline{T_{k}(\rho)}\|_{L^{2}((0,T)\times\R^{N})}\rightarrow_{k\rightarrow+\infty}0$$
For poving the previous inequality, we see that
$\|T_{k}(\rho)-\overline{T_{k}(\rho)}\|_{L^{1}((0,T)
\times\R^{N})}\rightarrow0$ for $k\rightarrow+\infty$. We then
conclude by interpolation with
$T_{k}(\rho)-\overline{T_{k}(\rho)}\in L^{q}(0,T)\times\R^{N})$ with
$q>2$. By H\"older inequality we obtain that:
$$\int^{t_{2}}_{t_{1}}
\int_{\R^{N}}(T_{k}(\rho)-\overline{T_{k}(\rho)}){\rm
div}u\,dxdt\rightarrow0\;\;\;\mbox{for}\;k\rightarrow+\infty.$$ We
have then:
\begin{equation}
\begin{aligned}
&\lim_{k\rightarrow+\infty}\int_{\R^{N}}\big(\overline{L_{k}(\rho)}-L_{k}(\rho)\big)(t_{2})dx-
\int_{\R^{N}}\big(\overline{L_{k}(\rho)}-L_{k}(\rho)\big)(t_{1})dx\\
&\hspace{3cm}=-\lim_{k\rightarrow+\infty}\int^{t_{2}}_{t_{1}}\int_{\R^{N}}\overline{T_{k}(\rho){\rm
div}u}-\overline{T_{k}(\rho)}{\rm div}u\,dxdt.
\end{aligned}
\end{equation}
We set:
$$
\begin{aligned}
&dft[\rho_{n}\rightarrow\rho](t)=\lim_{k\rightarrow+\infty}\int_{\R^{N}}\big(\overline{L_{k}(\rho)}-L_{k}(\rho)\big)(t)dx\\
&A(k,\rho)=\int^{t_{2}}_{t_{1}}\int_{\R^{N}}\overline{T_{k}(\rho){\rm
div}u}-\overline{T_{k}(\rho)}{\rm div}u\,dxdt.\\
\end{aligned}
$$
We can show as in the previous proof of theorem \ref{4principal}
that:
$$
\begin{aligned}
&\int^{t_{2}}_{t_{1}}\int_{K}\overline{T_{k}(\rho){\rm
div}u}-\overline{T_{k}(\rho)}{\rm div}u\,dxdt\\
&\hspace{2cm}=\lim_{n\rightarrow+\infty}\int^{t_{2}}_{t_{1}}
\int_{K}\big((P(\rho_{n})+\frac{\kappa}{2}\rho_{n}^{2})T_{k}(\rho_{n})-\overline{(P(\rho)+\frac{\kappa}{2}
\rho^{2})}\overline{T_{k}(\rho)}\big)dxdt.
\end{aligned}
$$
for any compact $K\subset\R^{N}$.
\\
Using the lemma \ref{4lemme2} we deduce that:
$$\lim_{n\rightarrow+\infty}\int^{t_{2}}_{t_{1}}\int_{\R^{N}}(P_{1}(\rho_{n})+\frac{\kappa}{2}
\rho_{n}^{2})T_{k}(\rho_{n})-(\overline{P_{1}(\rho)}+\frac{\kappa}{2}\overline{\rho^{2}})\overline{T_{k}(\rho)}dxdt)
\leq0.$$
We have then:
$$
\begin{aligned}
dft[\rho_{n}\rightarrow\rho](t_{2})-dft[\rho_{n}\rightarrow\rho](t_{1})\leq\lim_{k\rightarrow+\infty}
\biggl(\lim_{n\rightarrow+\infty}\int^{t_{2}}_{t_{1}}\int_{\R^{N}}P_{2}&(\rho_{n})T_{k}(\rho_{n})\\
&-\overline{P_{2}(\rho)}\,\overline{T_{k}(\rho)}dxdt\biggl).\\
\end{aligned}
$$
As the sequence $(\rho_{n})_{n\in\mathbb{N}}$ is bounded in
$L^{\infty}(0,T,L^{2}(\R^{N}))$, and $P_{2}$ is a bounded function,
we have:
$$
\begin{aligned}
&\lim_{k\rightarrow+\infty}\biggl(\lim_{n\rightarrow+\infty}\int^{t_{2}}_{\tau_{1}}\int_{\R^{N}}
P_{2}(\rho_{n})T_{k}(\rho_{n})-\overline{P_{2}(\rho)}\overline{T_{k}(\rho)}dxdt\biggl)\\
&\hspace{5cm}=\lim_{n\rightarrow+\infty}\int^{t_{2}}_{t_{1}}\int_{\R^{N}}P_{2}(\rho_{n})\rho_{n}
-\overline{P_{2}(\rho)}\,\overline{\rho}dxdt.
\end{aligned}
$$
Since the function $P_{2}$ is twice continuously differentiable and
compactly supported in $[0,+\infty)$, there exists $\Lambda>0$ big
enough such that both $\rho\rightarrow\Lambda\rho \log\rho-\rho
P_{2}(\rho)$ and $\rho\rightarrow\Lambda\rho \log\rho+\rho
P_{2}(\rho)$ are convex functions of $\rho$, indeed the twice
derivative are
positive.\\
As a consequence of weak lower semi-continuity of convex
functionals, we obtain:
$$
\begin{aligned}
&\lim_{n\rightarrow+\infty}\int^{\tau_{2}}_{t_{1}}\int_{\R^{N}}P_{2}(\rho_{n})\rho_{n}-\overline{P_{2}(\rho)}\,
\overline{\rho}\,dx\,dt\\
&\hspace{2cm}\leq\Lambda\int^{\tau_{2}}_{t_{1}}\int_{\R^{N}}(\overline{\rho\,\log\rho}-\rho\,\log\rho)dxdt+
\int^{\tau_{2}}_{t_{1}}\int_{\R^{N}}(P_{2}(\rho)-\overline{P_{2}(\rho)})\rho\, dx\,dt.\\
\end{aligned}
$$
Futhermore we have:
$$
\begin{aligned}
\int^{t_{2}}_{t_{1}}\int_{\R^{N}}(P_{2}(\rho)-\overline{P_{2}(\rho)})\rho dxdt&\leq\int^{\tau_{2}}_{\tau_{1}}\int_{\rho\leq\rho_{r}}(P_{2}(\rho)-\overline{P_{2}(\rho)})\rho\\
&\leq\Lambda\int^{t_{2}}_{t_{1}}\int_{\rho\leq\rho_{r}}(\overline{\rho\,\log\rho}-\rho\,\log\rho)dxdt\\
&
\leq\Lambda\rho_{r}\int^{t_{2}}_{t_{1}}\int_{\R^{N}}(\overline{\rho\,\log\rho}-\rho\,\log\rho)dxdt
\end{aligned}
$$
The previous relation gives:
$$dft[\rho_{n}\rightarrow\rho](t_{2})\leq dft[\rho_{n}\rightarrow\rho](t)+\omega\int^{t_{2}}
_{t_{1}}dft[\rho_{n}\rightarrow\rho](t).$$ Applying Gr\"onwall's
lemma we infer:
$$dft[\rho_{n}\rightarrow\rho](t_{2})\leq dft[\rho_{n}\rightarrow\rho](t_{1})exp(\omega(t_{2}-t_{1})).$$
We conclude that $dft[\rho_{n}\rightarrow\rho](t)=0\;\forall t$,
because $\rho^{0}_{n}$ converges strongly in $L^{1}$ to $\rho_{0}$.
\section{Weak solutions with data close to a stable equilibrium}
We consider in this section one situation which is rather different
from the three cases considered in the preceeding sections. This
situation is relevant for practical applications and realistic flow
and they involve conditions at infinity different from those
studied.\\
We wish to investigate the system $(NSK)$ with hypothesis close from
these of the theorem for strong solutions. We want then to study the
system with a density close from a stable equilibrium in the goal to
can choose initial data which avoid the vacuum.
We look now for a solution $(\rho,u)$ defined on $\R\times\R^{N}$ of
the system $(NSK)$ (where $P(\rho)=a\rho^{\gamma}$) with $\rho\geq0$ on $\R\times\R^{N}$.\\
In addition we require $(\rho,u)$ to satisfy the following limit
conditions:
$$(\rho,u)(x,t)\rightarrow
(\bar{\rho},0)\;\;\mbox{as}\;\;|x|\rightarrow +\infty,\;\;\mbox{for
all}\;\;t>0$$ where $\bar{\rho}>0$.
\\
Such an analysis requires the use of the Orlicz spaces. We define
the Orlicz space $L^{q}_{p}(\R^{N})$ as follows:
$$L^{q}_{p}(\R^{N})=\{f\in L^{1}_{loc}(\R^{N})/f 1_{\{|f|\leq\delta\}}\in L^{p}(\R^{N}),\;\;
f 1_{\{|f|\geq\delta\}}\in L^{q}(\R^{N})\}.$$ For the following
properties, we can give this definition too.
\begin{definition}
We define $\Psi$ as a convex function on $[0,+\infty)$ which is
equal (or equivalent) to $x^{p}$ for $x$ small and to $x^{q}$ for
$x$ large.
$$L^{q}_{p}(\R^{N})=\{f\in L^{1}_{loc}(\R^{N})/\Psi(f)\in L^{1}(\R^{N})\}.$$
We can check that $L^{q}_{p}(\R^{N})$ is a linear space. Now we
endow $L^{q}_{p}(\R^{N})$ with  a norm so that $L^{q}_{p}(\R^{N})$ is a separable Banach space.\\
$$\|f\|_{L^{q}_{p}(\R^{N})}=\inf\big(t>0/\;\;\Psi(\frac{f}{t})\leq 1\big).$$
\end{definition}
We recall now some useful properties of the Orlicz space:
\begin{corollaire}
We have:
\begin{enumerate}
\item Embedding:
$$L^{q}_{p}(\R^{N})\subset L^{q_{1}}_{p_{1}}(\R^{N})\;\;if\;\;1\leq
q_{1}\leq q<+\infty,\;\;1\leq p\leq p_{1}<+\infty.$$
\item Topology: $f_{n}\rightarrow_{n} 0$ in $L^{q}_{p}(\R^{N})$ if and only if
$\Psi(f_{n})\rightarrow_{n} 0$ and that:
$$\Psi(\frac{f}{\|f\|_{L^{q}_{p}(\R^{N})}})=1\;\;\;\mbox{if}\;\;f\ne
0$$
\end{enumerate}
\end{corollaire}
We recall now some useful properties of the Orlicz space:
\begin{proposition}
\label{4compositionorlicz} We then have:
\begin{itemize}
\item Dual space: If $p>1$ then
$(L^{q}_{p}(\R^{N}))^{'}=L^{q^{'}}_{p^{'}}(\R^{N})$ where
$q^{'}=\frac{q}{q-1},\,p^{'}=\frac{p}{p-1}.$
\item Let $F$ be a continuous function on $\R$ such that $F(0)=0$,
$F$ is differentiable at $0$ and
$F(t)|t|^{-\theta}\rightarrow\alpha\ne
0$ at $t\rightarrow +\infty$. Then if $q\geq\theta$,\\
$$F(f)\in L^{\frac{q}{\theta}}_{p}(\R^{N})\;\;\mbox{if}\;\;f\in
L^{q}_{p}(\R^{N}).$$
\end{itemize}
\end{proposition}
Our goal is now to get energy estimate. We have to face a new
difficulty. Indeed $\rho,\,\rho|u|^{2},\,\rho^{\gamma}$ need not
belong to
$L^{1}$.\\
We first want to explain how it is possible to obtain natural a
priori bounds which correspond to energy-like identities.\\
Next we write the following formal identities:
\begin{equation}
\begin{aligned}
&\frac{1}{\gamma-1}\frac{d}{dt}(\rho^{\gamma}-\bar{\rho}^{\gamma}
-\gamma\bar{\rho}^{\gamma-1}(\rho-\bar{\rho})+{\rm
div}[u\frac{\gamma}{\gamma-1}(\rho^{\gamma}-\bar{\rho}^{\gamma-1}\rho)]=u\cdot\n(\rho^{\gamma})\\[2mm]
&\rho\frac{d}{dt}\frac{|u|^{2}}{2}+\rho
u\cdot\n\frac{|u|^{2}}{2}-\mu\D u\cdot u-\xi\n{\rm
div}u\cdot u+au\cdot\n\rho^{\gamma}=\kappa\rho u\n(\phi*\rho-\rho).\\
\label{426}
\end{aligned}
\end{equation}
We may then integrate in space the equality (\ref{426}) and we get:
\begin{equation}
\begin{aligned}
&\big(\int_{\R^{N}}\rho\frac{|u|^{2}}{2}+\frac{a}{\gamma-1}(\rho^{\gamma}+(\gamma-1)\bar{\rho}^{\gamma}
-\gamma\bar{\rho}^{\gamma-1}\rho)+E_{global}[\rho-\bar{\rho}]dx\big)(t)\\
&+\int_{0}^{t}ds\int_{\R^{N}}2\mu|D u|^{2}+2\xi|{\rm
div}u|^{2}dx\leq\int_{\R^{N}}\rho_{0}\frac{|u_{0}|^{2}}{2}+\frac{a}{\gamma-1}(\rho_{0}^{\gamma}+(\gamma-1)\bar{\rho}^{\gamma}
-\gamma\bar{\rho}^{\gamma-1}\rho_{0})\\
&\hspace{11cm}+E_{global}[\rho_{0}-\bar{\rho}]dx.\\
\end{aligned}
\end{equation}
\begin{notation}
In the sequel we will note:
$$j_{\gamma}(\rho)=\rho^{\gamma}+(\gamma-1)\bar{\rho}^{\gamma}
-\gamma\bar{\rho}^{\gamma-1}\rho.$$
\end{notation}
Now we can recall a theorem (see \cite{4L2}) on the Orlicz space
concerning this quantity:
\begin{theorem}
$j_{\gamma}(\rho)\in L^{1}(\R^{N})$ if and only if
$\rho-\bar{\rho}\in L^{\gamma}_{2}.$
\end{theorem}
In following this theorem and our energy estimate we get that
$\rho-\bar{\rho}\in L^{\infty}(0,T;L^{\gamma}_{2}(\R^{N}))$ for
all $T\in\R$.\\
Moreover we have:
$$E_{global}[\rho(t,\cdot)-\bar{\rho}](x)=\frac{\kappa}{4}\big(\rho-\bar{\rho})^{2}+\phi*(\rho-\bar{\rho})^{2}-
2(\rho-\bar{\rho})\,(\phi*(\rho-\bar{\rho})).$$ Then in using the
fact that $\rho-\bar{\rho}\in
L^{\infty}(0,T;L^{\gamma}_{2}(\R^{N}))\hookrightarrow L^{2}(\R^{N})
+L^{\gamma}(\R^{N})$ and interpolation on $\n\phi$, we get that $\rho-\bar{\rho}\in L^{\infty}(0,T;L^{2}(\R^{N}))$.\\
We may now turn to our compactness result. First of all, we consider
sequences of solutions $(\rho_{n},u_{n})$ of the system
corresponding to some initial conditions
$(\rho^{0}_{n},u^{0}_{n})$.\\
In using the above energy inequalities, we assume that
$j_{\gamma}(\rho^{0}_{n})$, $E_{global}[\rho^{0}_{n}-\bar{\rho}]$
and $\rho^{0}_{n}|u^{0}_{n}|^{2}$ are bounded in
$L^{\infty}(L^{1}(\R^{N}))$ and that $\rho^{0}_{n}-\bar{\rho}$
converges weakly in
$L^{\gamma}_{2}(\R^{N})$ to some $\rho_{0}-\bar{\rho}$.\\
\\
We now assume that:
$$j_{\gamma}(\rho_{n}),\,E_{global}[\rho_{n}-\bar{\rho}],\,\rho_{n}|u_{n}-\bar{u}|^{2}\;\;\mbox{are
bounded in}\;\;L^{\infty}(0,T,L^{1}(\R^{N})),$$ Moreover we have for
all $,T\in(0,+\infty)$ and for all compact sets $K\subset\R^{N}$:
$$\rho_{n}-\bar{\rho}\in L^{\infty}(L^{2}(\R^{N})\;\;\;\mbox{and}\;\;\;\rho_{n}\;\;\mbox{is bounded in}
\;\;L^{q}((0,T)\times K),$$ for some $q>s$.
$$
\begin{aligned}
&Du_{n}\;\;\mbox{is bounded in}\;\;L^{2}(\R^{N}\times
(0,T))\;\;\;\;\;\;\;\;\;\;\;\;\;\;\;\;\;\;\;\;\;\;\;\;\;\;\;\;\;\;\;\;\;\;\;\;\;\;\;\;\;\;\\
&u_{n}\;\;\mbox{is bounded
in}\;\;L^{2}(0,T,H^{1}(B_{R}))\;\;\mbox{for all}\;\;R,T\in
(0,+\infty)\;.\;\;\;\;\;\;\;\;\;\;\;\;\;\;\;\;\;\;\;\;\;\;\;\;\;\;\;\;\;\;\;\;\;\;\;\;\;\;\;\;\;\\
\end{aligned}
$$
\\
Extracting subsequences if necessary, we may assume that
$\rho_{n},\,u_{n}$ converge weakly respectively in
$L^{2}((0,T)\times B_{R})$, $L^{2}(0,T;H^{1}(B_{R}))$ to $\rho,\,u$
for all $R,T\in(0,+\infty)$. We also extract subsequences for which
$\rho_{n}u_{n},\,\rho_{n}u_{n}\otimes u_{n}$ converge weakly as
previously.
\begin{rem}
We notice that the conditions at infinity are implicitly contained
in the fact that $(\rho_{n}-\bar{\rho})^{2}$ and
$\rho_{n}|u_{n}|^{2}\in L^{1}(\R^{N})$.
\end{rem}
We then have the following theorem.
\begin{theorem}
Let $\gamma\geq1$. We assume that $\rho^{0}_{n}$ converges in
$L^{1}(B_{R})$ (for all $R\in(0,+\infty)$) to $\rho_{0}$. Then
$(\rho_{n},u_{n})_{n\in\mathbb{N}}$ converges in distribution sense
to $(\rho,u)$
a solution of $(NSK)$.\\
Moreover we have for all $R,T\in(0,+\infty)$:
$$
\begin{aligned}
&\rho_{n}\rightarrow_{n}\rho\;\;\mbox{in}\;\;C([0,T],L^{p}(B_{R}\times(0,T))\cap
L^{s_{1}}(B_{R}\times(0,T))\:\;\mbox{for all}\;\;1\leq p<s,1\leq
s_{1}<q.\\
\end{aligned}
$$
with $q=s+\frac{4}{N}-1$.
\end{theorem}
{\bf Proof:}
\\
As in the theorem \ref{4principal}, we want test the strong
convergence of $\rho_{n}$ on concav function $B$. Since the proof is
purely local, we have again for small enough $\e>0$:
\begin{equation}
\begin{aligned}
&(\rho_{n})^{\e}\big((\mu+\xi){\rm
div}u_{n}-a(\rho_{n})^{\gamma}-\frac{\kappa}{2}\rho_{n}^{2}\big)\rightharpoonup_{n}\overline{(\rho)^{\e}}
\big((\mu+\xi){\rm
div}u-a\overline{\rho^{\gamma}}-\frac{\kappa}{2}\overline{\rho^{2}}\big)\;\;\;\\
&\hspace{11cm}\mbox{in}\;\;{\cal D}^{'}\big((0,\infty)\times\R^{N}\big),\\[2mm]
&\frac{d}{dt}(\overline{\rho^{\e}})+{\rm
div}(u\overline{\rho^{\e}})\geq (1-\e)({\rm
div}u)\overline{\rho^{\e}}\;\;\;\mbox{in}\;\;{\cal D}^{'}\big((0,\infty)\times\R^{N}\big).\\
\end{aligned}
\label{4H40}
\end{equation}
Next since $(\overline{\rho^{\e}})^{\frac{1}{\e}}\in
L^{2}(B_{R}\times(0,T))$ for all $R,T\in(0,+\infty)$, as in the
theorem \ref{4principal} in using a result of type Diperna-Lions on
renormalized solutions, we get:
\begin{equation}
\label{4H41}
\frac{d}{dt}(\overline{(\rho)^{\e}}^{\frac{1}{\e}})+{\rm
div}(u\overline{(\rho)^{\e}}^{\frac{1}{\e}})\geq
0\;\;\;\mbox{in}\;\; {\cal D}^{'}\big((0,\infty)\times\R^{N}\big).
\end{equation}
while we have $\overline{(\rho)^{\e}}^{\frac{1}{\e}}\leq\rho$ a.e in
$\R^{N}\times(0,+\infty)$ and
$\overline{(\rho)^{\e}}^{\frac{1}{\e}}_{/t=0}=\rho_{/t=0}$ in
$\R^{N}$.\\
Now in subtracting the second equality of (\ref{4H41}) from the mass
equation and setting $f=\rho-\overline{\rho^{\e}}^{\frac{1}{\e}}$,
we have:
\begin{equation}
\frac{d}{dt}(f)+{\rm
div}(uf)\leq0,\;f\geq0\;\;\mbox{a.e},\;f_{/t=0}=0\;\;\mbox{in}\;\;
\R^{N}. \label{4H42}
\end{equation}
Next we want again to show from (\ref{4H42}) that $f=0$, in
integrating (\ref{4H42}) and in using the fact that $f\leq0$ to
conclude that $f=0$. The difference with the proof of theorem
\ref{4principal}
is to justify the integration by parts as we work in different energy space.\\
We need a cut-off function. We introduce $\varphi\in
C^{\infty}_{0}(\R^{N}),\,0\leq\varphi\leq1,\,\mbox{supp}\varphi\subset
B_{2} ,\,\varphi=1\;\mbox{on}\;B_{1}$ and we set
$\varphi_{R}=\varphi(\frac{x}{R})$ for $R\geq1$. Multiplying
(\ref{4H42}) by $\varphi_{R}(x)$, we obtain:
\begin{equation}
\begin{aligned}
&\frac{d}{dt}\int_{\R^{N}}f\varphi_{R}(x)dx=&\int_{\R^{N}}\frac{1}{R}fu\cdot\n\varphi(\frac{x}{R}).\\
\end{aligned}
\label{4H43}
\end{equation}
Next, if $T>0$ is fixed, we see that
$\mbox{supp}\,\n\varphi(\frac{\cdot}{R})\subset\{R\leq|x|\leq2R\}$ ,
therefore, for $R$ large enough we have:
\begin{equation}
\begin{aligned}
\frac{d}{dt}\int_{\R^{N}}f\varphi_{R}(x)dx=&\int_{\R^{N}}\frac{1}{R}fu.\n\varphi(\frac{x}{R}),\\
&\leq\frac{C}{R}\int_{\R^{N}}f|u|1_{(R\leq|x|\leq2R)}dx,\;\;\mbox{for}\;t\in(0,T),\\
\end{aligned}
\label{4H44}
\end{equation}
To conclude that $f=0$, we only have to prove that:
\begin{equation}
\frac{1}{R}\int_{\R^{N}}f|u|1_{(R\leq|x|\leq2R)}dx\rightarrow0\;\;\mbox{as}\;R\rightarrow+\infty.
\label{4H45}
\end{equation}
We now use that fact that $f\in L^{\infty}(0,T,L^{2}(\R^{N}))$ and
$f|u|^{2}\in L^{\infty}(0,T,L^{1}(\R^{N})$ for all $T\in(0,+\infty)$
to control (\ref{4H45}). The second fact is obvious since
$0\leq\rho$ and
$\rho|u|^{2}\in L^{\infty}(0,T;L^{1}(\R^{N}))$.\\
In order to prove the first claim, we only have to show that
$\overline{(\rho)^{\e}}^{\frac{1}{\e}}-\bar{\rho}\in
L^{\infty}(0,T,L^{2}(\R^{N}))$. We rewrite
$(\rho_{n})^{\e}-(\bar{\rho})^{\e}=(\bar{\rho}+(\rho_{n}-\bar{\rho}))^{\e}-
(\bar{\rho})^{\e}$ is bounded in $L^{\infty}(0,T,L^{2}(\R^{N}))$ in
using proposition \ref{4compositionorlicz} with
$F(x)=(\bar{\rho}+x)^{\e}-(\bar{\rho})^{\e}$.
So we have $\sqrt{f}\in L^{\infty}(L^{4}(\R^{N}))$ and we get:
$$
\begin{aligned}
\frac{1}{R}\int^{T}_{0}dt\int_{\R^{N}}f|u|&1_{(R\leq|x|
\leq2R)}dx\leq \frac{C_{0}}{R}\mbox{meas}
(C(0,R,2R))^{\frac{1}{4}}.\\
\end{aligned}
$$
We recall that:
$$\mbox{meas}(C(0,R,2R))\thicksim_{R\rightarrow+\infty}C(n)R^{N}$$
Then we get:
$$\frac{d}{dt}\int_{\R^{N}}f\varphi_{R}(x)dx\rightarrow_{R\rightarrow
+\infty}0$$ and we conclude as in the proof of theorem
\ref{4principal}.\\
At this stage, it only remains to show that, for instance,
$\rho_{n}$ converge to $\rho$ in $C([0,T],L^{1}(B_{R}))$ for all
$R,T\in(0,+\infty)$. In order to do so, we just have to localize the
corresponding argument in the proof of theorem \ref{4principal}.\\
Therefore we choose for $R,T\in (0,+\infty)$ fixed, $\varphi\in
C_{0}^{\infty}(\R^{N})$ such that $\varphi=1$ on $B_{R}$,
$0\leq\varphi$ on $\R^{N}$. Then, we observe that we have:
$$
\begin{aligned}
&\frac{\p}{\p t}(\va^{2}\rho_{n})+{\rm
div}(u_{n}(\va^{2}\rho_{n}))=\rho_{n}u_{n}\cdot\n\va^{2},\;\frac{\p}{\p
t}(\va^{2}\rho)+{\rm
div}(u(\va^{2}\rho))=\rho u\cdot\n\va^{2}\\
&\frac{\p}{\p t}(\va\sqrt{\rho_{n}})+{\rm
div}(u_{n}(\va\rho_{n}))=\frac{1}{2}({\rm
div}u_{n})\va\sqrt{\rho_{n}}+\sqrt{\rho_{n}}u_{n}\cdot\n\va,\\
&\frac{\p}{\p t}(\va\sqrt{\rho})+{\rm
div}(u(\va\rho))=\frac{1}{2}({\rm
div}u)\va\sqrt{\rho}+\sqrt{\rho}u\cdot\n\va.\\
\end{aligned}
$$
From these equations, we deduce as in the proof of the previous
theorem, that $\va^{2}\rho\in C([0,+\infty),L^{1}(\R^{N}))$,
$\va\sqrt{\rho}\in  C([0,+\infty),L^{2}(\R^{N}))$ and that
$\va\sqrt{\rho_{n}}$ converges weakly in $L^{2}(\R^{N})$, uniformly
in $t\in[0,T]$. Therefore, in order to conclude, we just have to
show that we have:\\
$$\int_{\R^{N}}\va^{2}\rho_{n}(t_{n})dx\rightarrow\int_{\R^{N}}\va^{2}\rho(\bar{t})dx$$
whenever $t_{n}\in[0,T]$, $t_{n}\rightarrow_{n}\bar{t}$, and this is
straightforward since we have, in view of the above equation:\\
$$
\begin{aligned}
\int_{\R^{N}}\va^{2}\rho_{n}(t_{n})dx&=\int_{\R^{N}}\va^{2}(\rho_{0})_{n}dx+\int^{t_{n}}_{0}ds\int_{\R^{N}}
\rho_{n}u_{n}\cdot\n\va^{2}\,dx\\
&\longrightarrow_{n}\int_{\R^{N}}\va^{2}\rho_{0}dx+\int^{T}_{0}ds\int_{\R^{N}}
\rho u\cdot\n\va^{2}\,dx=\int_{\R^{N}}\va\rho(\bar{t})dx.\\
\end{aligned}
$$
\hfill{$\Box$}


\begin{thebibliography}{}
\bibitem{4NS}
A. Novotn\'y and I. Stra$\check{\mbox{s}}$kraba. Introduction to the
mathematical theory of compressible flow,
Oxford lecture series in mathematics and its application. 27\\
\end{thebibliography}
\end{document}